\documentclass[a4paper]{article}
\usepackage{hyperref}
\usepackage{amsmath}
\usepackage{amsfonts}
\usepackage{amsthm}
\usepackage{indentfirst}
\usepackage{geometry}
\usepackage{color}
\newtheorem{prop}{Proposition}[section]
\renewcommand{\theequation}{\thesection.\arabic{equation}}
\usepackage{indentfirst}
\date{}
\makeatletter \@addtoreset{equation}{section}

\newtheorem{theorem}{Theorem}[section]
\newtheorem{lemma}[theorem]{Lemma}

\newtheorem{rem}{Remark}[section]
\newtheorem{de}{Definition}[section]

\begin{document}
	
	\begin{center}
		{\Large
			Sharp criteria for a degenerate diffusion-aggregation system with the intermediate exponent}
	\end{center}
	
	\vskip 5mm
	
	\begin{center}
		{\sc Tiantian Zhou \quad Li Chen \quad Yutian Lei}
	\end{center}

	\vskip 5mm {\leftskip5mm\rightskip5mm \normalsize
		\noindent{\bf{Abstract}}
		In this paper, we investigate a multi-dimensional nonlocal degenerate  diffusion-aggregation equation with a diffusion exponent $m$ in the intermediate range $\frac{2d}{2d-\gamma}<m<\frac{d+\gamma}{d}$, where the nonlocal aggregation term is given by singular potential $|x|^{-\gamma}$, $0<\gamma\leq d-2$. Under two different assumptions on the initial data, we establish two sharp criteria (i.e., the critical thresholds in Theorem 1.1 and Theorem 1.2) governing the global existence and finite-time blow-up of solutions. Once the initial free energy is less than a constant that depends on the total mass (or depends on the extremum function of the Hardy-Littlewood-Sobolev inequality), the first criterion depends on the relationship between the $L^{\frac{2d}{2d-\gamma}}$-norm of initial data and total mass, while the second relies on the relationship between the $L^m$-norm of initial data and extremal function.
		In the discussion of the second criterion, we do not require $L^\infty(\mathbb{R}^d)$ boundedness of the initial data, which is necessary in reference \cite{B}. Furthermore, with the help of moment estimate, we manage to prove the compactness argument on the whole space by using the Lions-Aubin Lemma.
		Importantly, we demonstrate that the two initial free energy conditions on which two criteria are based are equivalent. Building on this, we further prove that the two sharp criteria themselves  are also equivalent, thereby unifying the classification results obtained from two different approaches.
		
		}

	\renewcommand{\theequation}{\thesection.\arabic{equation}}
	\catcode`@=11
	\@addtoreset{equation}{section}
	\catcode`@=12
	
	\section{Introduction}
	In this paper, we study the following degenerate diffusion-aggregation nonlocal equation in spatial dimension $d\geq 3$
	\begin{equation}\label{a1}
		\begin{cases}
			\rho_t=\Delta\rho^m-\text{div}(\rho\nabla c), &x\in \mathbb{R}^d, t\geq 0,\\
			c=\bar{c}_d\frac{1}{|x|^{\gamma}}\ast  \rho, &x\in \mathbb{R}^d, t\geq 0,\\
			\rho(x,0)=\rho_0(x), &x\in \mathbb{R}^d.
		\end{cases}
	\end{equation}
	where the diffusion exponent is taken to be $m\in (m_*,m^*)$ for $m_*=\frac{2d}{2d-\gamma}$, $m^*=\frac{d+\gamma}{d}$ and $0<\gamma \leq d-2$,  $\bar{c}_d=\frac{1}{\gamma d \alpha(d)}$ and $\alpha(d)=\frac{\pi^{d/2}}{\Gamma(1+d/2)}$ is the volume of the unit ball of $\mathbb{R}^d$. Here $\rho$ represents the bacteria density, and $c$ represents the chemical substance concentration. This model is commonly applied to characterize the collective movement of cells or the evolution of bacterial density.
	
	The system \eqref{a1} has two fundamental invariances:
	One is that the total mass of cell always remains constant for all times. Namely,
	\begin{equation}\label{mass}
		\int_{\mathbb{R}^d}\rho(x,t)dx=\int_{\mathbb{R}^d}\rho_0(x)dx \quad \text{for } t\geq 0.
	\end{equation}
	A critical diffusion exponent $m^*$ comes from the scaling invariance of the total mass. More precisely, if $\rho$ is a solution of \eqref{a1}, then the mass invariant function $\rho_\lambda(x,t):=\lambda^d\rho(\lambda x, \lambda^\nu t)$ is also a solution of \eqref{a1} if and only if $m=m^*$ and $\nu=\gamma+2$.
	
	Another important quantity is the free energy, which is defined as
	\begin{equation}\label{F}
		\mathcal{F}(\rho) = \frac{1}{m-1} \int_{\mathbb{R}^d} \rho^m(x,t) dx - \frac{1}{2} \int_{\mathbb{R}^d} \rho(x,t) c(x,t) dx.
	\end{equation}
	The free energy decays in time according to the entropy-entropy production relation:
	$$
	\frac{d}{dt} \mathcal{F}(\rho(\cdot,t)) + \int_{\mathbb{R}^d} \rho \left| \nabla \left( \frac{m}{m-1} \rho^{m-1} - c \right) \right|^2 dx = 0.
	$$
	The other critical exponent $m_*$ comes from the conformal invariance properties of the free energy $\mathcal{F}$. More precisely, if  $\rho$ is a solution of \eqref{a1}, then the conformal invariant function $\rho_\lambda(x,t):=\lambda^{\frac{d}{m}}\rho(\lambda x, \lambda^\nu t)$ is also a solution of \eqref{a1} if and only if $m=m_*$ and $\nu=\frac{\gamma}{2}+2$.
	
	When $\gamma=d-2$, the convolution of the Riesz potential reduces to that of the Newtonian potential. Therefore, the system \eqref{a1} reduces to the following well-studied degenerate Patlak-Keller-Segel system
	\begin{equation}\label{a3}
		\begin{cases}
			\rho_t=\Delta\rho^m-\text{div}(\rho\nabla c), &x\in \mathbb{R}^d, t\geq 0,\\
			-\Delta c= \rho, &x\in \mathbb{R}^d, t\geq 0,\\
			\rho(x,0)=\rho_0(x), &x\in \mathbb{R}^d.
		\end{cases}
	\end{equation}
	
	As the diffusion exponent $m=1$ and the space dimension $d=2$ in \eqref{a3}, the above system, which has been extensively investigated in \cite{BDP,PB}, corresponds to the Patlak-Keller-Segel system or to the classical Smoluchowski-Poisson system in two dimensions with linear diffusion \cite{K,P}. Actually, Dolbeault and Perthame \cite{DP} established that the qualitative behavior of solutions depends solely on the initial mass, i.e., the solution exists globally if the total mass is less than $8\pi$, whereas the solution exhibits blow-up in finite time if the total mass is bigger than $8\pi$. The study of the critical mass case $8\pi$ in \cite{BCM} revealed that while the solution exists globally, it blows up by concentrating into a Dirac measure at the center of mass as $t\to \infty$.
	This phenomenon holds considerable biological and physical relevance, as it is directly connected to the emergence of cellular aggregates in biology and to gravitational collapse in astrophysics.
	
	For dimension $d\geq 3$ in \eqref{a3}, increasing attention has been directed toward the sharp threshold behavior in the degenerate parabolic-elliptic Keller-Segel system. A critical diffusion exponent $m^*=2-\frac{2}{d}$ in \eqref{a3} arises from the scaling invariance of the total mass, as identified by Blanchet, Carrillo and Laurencot \cite{BCL}, where Sugiyama \cite{S} and Sugiyama and Kunii \cite{SK} classified the solutions based on this exponent: if $m>m^*$, the solution exists globally for any initial data and if $1<m\leq m^*$, the global existence or finite-time blow-up of solutions is determined by whether the size of the initial data exceeds the critical threshold. For more results on this exponent, we refer to \cite{BL, Blan, LS, S07} and the references therein. Another critical exponent $m_*=\frac{2d}{d+2}$ in \eqref{a3}, which was discovered by Chen, Liu and Wang \cite{CLW}, comes from the conformal invariance properties of the free energy. Their work shows that the $L^{m_*}$-norm of positive stationary solutions serves as a critical criterion for distinguishing between global existence and finite-time blow-up. Inspired by Chen and Wang \cite{CW}, we note that they utilized estimates of the solution in the $L^{m_*}$-norm, which are intimately connected to the optimal constant in the Hardy-Littlewood-Sobolev inequality, to establish an explicit criterion for the global existence and finite-time blow-up of solutions to \eqref{a3} in the range $m_*<m<m^*$.
	
	The study of general potential functions has been built upon seminal contributions in the literature. Bedrossian \cite{B11} investigated the asymptotic behavior of global weak solutions under the condition that the potential satisfies $\gamma\in [d-2,d-1]$. In particular, Wang, Chen and Hong \cite{WCH} established the $L^m$-norm of the initial data as a critical criterion for the behavior of solutions to \eqref{a1} in the conformally invariant case $m=m_*$. Specifically, solutions exist globally if this norm is below a universal threshold. Conversely, finite-time blow-up occurs when the $L^m$-norm of the initial data exceeds that of the stationary solution and the initial free energy is below that of the stationary solution. This naturally leads to the central question of how to classify solution behaviors for the intermediate exponent range $m_*<m<m^*$. Therefore, the main objective of this work is to classify the global existence and finite-time blow-up of solutions based on their initial data, using two different methods. Both methods rely on the well-known Hardy-Littlewood-Sobolev (HLS) inequality \cite{LL}, which states as follows.
	
	\begin{prop}[HLS inequality]
		Let $\alpha_1, \alpha_2>1$ such that $\frac{1}{\alpha_1}+\frac{1}{\alpha_2}+\frac{\gamma}{d}=2$. Then
		\begin{equation}\label{hls}
			\iint_{\mathbb{R}^d \times \mathbb{R}^d} \frac{g(x) h(y)}{|x-y|^{\gamma}} dxdy \leq C_{\text{HLS}} \|g\|_{L^{\alpha_1}(\mathbb{R}^d)}\|h\|_{L^{\alpha_2}(\mathbb{R}^d)}
		\end{equation}
		for any $g \in L^{\alpha_1}(\mathbb{R}^d)$ and $h \in L^{\alpha_2}(\mathbb{R}^d)$, where  if $\alpha_1=\alpha_2=\frac{2d}{2d-\gamma}$, the above inequality becomes
		\begin{equation}\label{hls2}
			\iint_{\mathbb{R}^d \times \mathbb{R}^d} \frac{g(x) g(y)}{|x-y|^{\gamma}} dxdy \leq C_{\text{HLS}} \|g\|_{L^{\frac{2d}{2d-\gamma}}(\mathbb{R}^d)}^2,
		\end{equation}
		and the best constant is
		$$
		C_{\text{HLS}} = \pi^{\gamma/2} \frac{\Gamma(d/2-\gamma/2)}{\Gamma(d-\gamma/2)} \left\{ \frac{\Gamma(d/2)}{\Gamma(d)} \right\}^{-1+\gamma/d}.
		$$
		Moreover, there is equality in \eqref{hls} if and only if $g(x)= AU_{\lambda,\tilde{x}}$ and $h(x)=cg(x)$, for some constant $A>0$ and parameters $\lambda\in \mathbb{R}$, $\tilde x \in \mathbb{R}^d$, where
		$$
		U_{\lambda,\tilde x} = \left(\lambda^2+|x-\tilde x|^2 \right)^{(\gamma-2d)/2}.
		$$
		And
		$\|U_{\lambda, \tilde x}\|_{L^{m_*}(\mathbb{R}^d)}$ is a constant independent of $\lambda$ and $\tilde x$.
	\end{prop}
	Besides this inequality, let us recall the Gagliardo-Nirenberg-Sobolev inequality \cite{PB}
	\begin{equation}\label{GNS}
		\|\rho\|_{L^{q_1}(\mathbb{R}^d)}\leq S_d^{-\beta/2}\|\nabla \rho\|_{L^2(\mathbb{R}^d)}^{\beta}\|\rho\|_{L^{q_2}(\mathbb{R}^d)}^{1-\beta}.
	\end{equation}
	Here $\beta=\frac{2d(q_1-q_2)}{q_1[2d-(d-2)q_2]}$ with $q_2<q_1<\frac{2d}{d-2}$ and
	$S_d$ is given by the Sobolev inequality \cite{LL} on that
	\begin{equation}\label{So}
		S_d\|\rho\|_{L^{\frac{2d}{d-2}}(\mathbb{R}^d)}^2 \leq \|\nabla \rho\|_{L^2(\mathbb{R}^d)}^2, \quad S_d=\frac{d(d-2)}{4}2^{\frac{2}{d}}\pi^{1+\frac{1}{d}}\Gamma\left(\frac{d+1}{2}\right)^{-\frac{2}{d}}.
	\end{equation}

	The definition of weak solution is given in the following.
	\begin{de}\label{de}
		(Weak solution) Let $\rho_0 \in L^1_+(\mathbb{R}^d)$ be the given initial data.
		$\rho$ is a weak solution to the system \eqref{a1} if it satisfies
		$$
		\begin{aligned}
			&\int_{\mathbb{R}^d} \psi \rho(\cdot, t) dx-\int_{\mathbb{R}^d} \psi \rho_0(x) dx\\
			=&\int_0^t \int_{\mathbb{R}^d} \Delta \rho^m \psi dx ds
			-\frac{\bar{c}_d\gamma}{2} \int_0^t \iint_{\mathbb{R}^d \times \mathbb{R}^d} \frac{[\nabla \psi(x) - \nabla \psi(y)] \cdot (x - y)}{|x - y|^2} \frac{\rho(x, s) \rho(y, s)}{|x - y|^{\gamma}} dx dy ds
		\end{aligned}
		$$
		for any $\psi \in C_0^\infty(\mathbb{R}^d)$ and $0<t<\infty$.
	\end{de}
	
	In this paper, we focus on finding criteria on initial data that distinguish between global existence and finite-time blow-up of solutions.

	The first main result is given by the following theorem.
	\begin{theorem}\label{th1}
		Let $m\in (m_*,m^*)$ and $0<\gamma \leq d-2$. Assume the initial data satisfies
		\begin{equation}\label{a2}
			\rho_0\in L_+^1(\mathbb{R}^d)\cap L^m(\mathbb{R}^d),\quad  \nabla \rho_0^m\in L^2(\mathbb{R}^d),
		\end{equation}
		and initial free energy satisfies
		\begin{equation}\label{c1}
			\mathcal{F}(\rho_0)<\frac{\gamma-d(m-1)}{\gamma(m-1)} \left(\frac{2d^2\alpha(d)}{ C_{\text{HLS}}} \right)^{\frac{d(m-1)}{\gamma-d(m-1)}} \|\rho_0\|_{L^1(\mathbb{R}^d)}^{\frac{2d-m(2d-\gamma)}{\gamma-d(m-1)}}.
		\end{equation}
		The following holds
		
		(1) If $\|\rho_0\|_{L^{m_*}(\mathbb{R}^d)}<s_*$, then the weak solution $\rho$ exists globally with
		$$
		\rho \in L^\infty \left(0,\infty; L^1_+ \cap L^m(\mathbb{R}^d) \right),
		$$
		and for $T>0$,
		$$
		\nabla \rho \in L^2 \left(0,T; L^{\kappa}(\mathbb{R}^d) \right), \quad \kappa=\min \left\{2, \frac{2m}{3-m} \right\},
		$$
		$$
		\partial_t \rho \in L^2\left(0,T; W_{loc}^{-1,s}(\mathbb{R}^d)\right), \quad s\in \min \left\{\frac{2m}{m+1}, \frac{m(m+1)d}{(2d-\gamma-2)m+[m(\gamma+1-d)+d](m+1)}\right\}.
		$$

		(2) If $\|\rho_0\|_{L^{m_*}(\mathbb{R}^d)}>s_*$ and $\int_{\mathbb{R}^d}|x|^2\rho_0(x)dx<\infty$, then the weak solution $\rho$ does not exist globally in time. Namely, there exists a $T^*>0$ such that it blows up in finite time $T^*$ in the following sense
		\begin{equation}\label{r55}
			\limsup_{t \to T^*} \|\rho(x, \cdot)\|_{L^m(\mathbb{R}^d)} = +\infty.
		\end{equation}
		Here $s_*$ is universal constants given by
		\begin{equation}\label{s}
			s_*:=\left(\frac{2d^2 \alpha(d)}{C_{\text{HLS}}}\right)^{\frac{\gamma}{2\gamma-2d(m-1)}}
			\|\rho_0\|_{L^1(\mathbb{R}^d)}^{\frac{\gamma m-2d(m-1)}{2\gamma-2d(m-1)}}.
		\end{equation}
	\end{theorem}
	
	\begin{rem}
		The case $\|\rho_0\|_{L^{m_*}(\mathbb{R}^d)}=s_*$ does not happen if \eqref{c1} holds. Hence, the above classification yields a sharp criterion in the sense that whether the $L^{m_*}$-norm of the initial data satisfying \eqref{c1} exceeds the threshold $s_*$ strictly determines global existence or finite-time blow-up of the solutions.
	\end{rem}

	Another classification criterion, based on the extremal function of the HLS inequality, has also attracted considerable attention. Ogawa \cite{O} established that the classification criterion for solutions to \eqref{a3} at the critical exponent $m = m_*$ is fundamentally rooted in the extremal function of the Sobolev inequality; through duality, this criterion equivalently relies on the extremal function of the HLS inequality. For the intermediate exponent range $m \in (m_*, m^*)$, Kimijima, Nakagawa, and Ogawa \cite{KNO} obtained a classification of global existence and finite-time blow-up for weak solutions to \eqref{a3} based on scale-invariant norms of the initial data. They showed that the critical threshold governing the solution behavior is determined by the extremal function associated with the optimal constant in the modified HLS inequality. Subsequently, Bian \cite{B} established the $L^m$-norm of the stationary solutions, which coincides with this extremal function, as a criterion for the behavior of solutions to \eqref{a1} with $m=m_*$. More recently, Bian and Du \cite{BD} applied this approach to classify solutions of \eqref{a1} for intermediate exponents, where Moser iteration \cite{BLZ, CW19, LW} plays a crucial role in deriving $L^\infty$ estimate and compactness argument.
	
	We investigate the extremal function for attaining the optimal constant in a modified HLS inequality, with the aim of using its associated invariant norm as a sharp criterion to classify initial data as leading to either global existence or blow-up.
	
	For $1<\frac{2d}{2d-\gamma}<m$, it holds from the interpolation inequality
	\begin{equation}\label{int}
		\|\rho\|_{L^{m_*}(\mathbb{R}^d)} \leq \|\rho\|_{L^1(\mathbb{R}^d)}^{1-\theta} \|\rho\|_{L^m(\mathbb{R}^d)}^{\theta}, \quad \theta=\frac{m\gamma}{2d(m-1)}.
	\end{equation}
	Thus, \eqref{hls2} yields that
	\begin{equation}\label{bhls}
		\int_{\mathbb{R}^d}\int_{\mathbb{R}^d}\frac{\rho(x,t)\rho(y,t)}{|x-y|^{\gamma}}dxdy\leq C_{\text{HLS}}\|\rho\|_{L^1(\mathbb{R}^d)}^{2(1-\theta)}\|\rho\|_{L^m(\mathbb{R}^d)}^{2\theta},
	\end{equation}
	and
	\begin{equation}\label{b1}
		C_{\text{HLS}} := \sup\Bigg\{\frac{\int_{\mathbb{R}^d}\int_{\mathbb{R}^d}\frac{\rho(x,t)\rho(y,t)}{|x-y|^{\gamma}}dxdy}{\|\rho\|_{L^1(\mathbb{R}^d)}^{2(1-\theta)} \|\rho\|_{L^m(\mathbb{R}^d)}^{2\theta}}\; \bigg|\; \rho \in L^1(\mathbb{R}^d) \cap L^{m}(\mathbb{R}^d), \rho\neq 0\Bigg\}<\infty.
	\end{equation}
	According to \cite[Proposition 4.1]{KNO} (see also \cite{BD}), we know that there exists the  function $V\in L^1(\mathbb{R}^d) \cap L^{m}(\mathbb{R}^d)$
	that attains the maximizing problem \eqref{b1} and it satisfies that
	\begin{equation}\label{b2}
		C_{\text{HLS}}\|V\|_{L^1(\mathbb{R}^d)}^{2(1-\theta)}\|V\|_{L^m(\mathbb{R}^d)}^{2\theta}=\int_{\mathbb{R}^d}\int_{\mathbb{R}^d}\frac{V(x)V(y)}{|x-y|^{\gamma}}dxdy.
	\end{equation}

	Based on the above preparation, our second result of this paper is as follows:
	\begin{theorem}\label{th2}
		Let $m \in (m_*, m^*)$ and
		\begin{equation}\label{gamma}
			\begin{cases}
				0<\gamma \leq d-2, & d=3,4,\\
				\frac{d}{2}-2\leq \gamma\leq d-2, &d\geq 5.
			\end{cases}
		\end{equation}
		Assume the initial data  satisfies
		\begin{equation}\label{in}
			\rho_0\in L_+^1(\mathbb{R}^d)\cap L^{p_0}(\mathbb{R}^d), \quad \int_{\mathbb{R}^d}|x|^2\rho_0 dx<\infty, \quad
			\nabla \rho_0^m\in L^2(\mathbb{R}^d),
		\end{equation}
		where $p_0=\max\{2m+2,\gamma+2, \frac{2md}{d-m(d-\gamma-2)}-1+\tau, d-\frac{dm}{2}-1+\tau\}$ for any $0<\tau \ll  1$, and the initial free energy satisfies
		\begin{equation}\label{c2}
			\|\rho_0\|_{L^1(\mathbb{R}^d)}^\frac{\tau_1 m}{\tau_2} \mathcal{F}(\rho_0) < \|V\|_{L^1(\mathbb{R}^d)}^{\frac{\tau_1m}{\tau_2}}\mathcal{F}(V),
		\end{equation}
		where $\tau_1$ and $\tau_2$ are given by
		\begin{equation}\label{tau}
			\tau_1=\frac{d-\gamma}{2-m}-\frac{d}{m}, \quad \tau_2=d-\frac{d-\gamma}{2-m}.
		\end{equation}
		Then the following statements hold
		
		(1) If
		$
		\|\rho_0\|_{L^m(\mathbb{R}^d)}<s^*,
		$
		then there exists a global weak solution $\rho$ to \eqref{a1}, which satisfies
		$$
		\rho\in L^\infty\left(0,\infty;L^{p_0}(\mathbb{R}^d)\right).
		$$
		
		(2) If
		$ \|\rho_0\|_{L^m(\mathbb{R}^d)}> s^*$,
		then the  weak solution of system \eqref{a1} blows up in finite time $T^*$ in the sense of \eqref{r55}.
		Here $s^*$ is universal constants given by
		\begin{equation}\label{ss}
			s^* =  \left(\frac{2d^2\alpha(d)}{ C_{\text{HLS}}}\right)^{\frac{d(m-1)}{\gamma m-md(m-1)}}\|\rho_0\|_{L^1(\mathbb{R}^d)}^{\frac{\gamma m-2d(m-1)}{\gamma m-md(m-1)}}.
		\end{equation}
	\end{theorem}

	\begin{rem}\label{re1}
		One can easily check that the case $\|\rho_0\|_{L^m(\mathbb{R}^d)}= s^*$ does not happen if \eqref{c2} holds. Hence, the above classification yields a sharp criterion in the sense that whether the $L^{m}$-norm of the initial data satisfying \eqref{c2} exceeds the threshold $s^*$ strictly determines global existence or finite-time blow-up of the solutions.
	\end{rem}
	
	\begin{rem}\label{re2}
		A key step in our proof is establishing the finite-time boundedness of the first-order moment of $\rho_\varepsilon^m$. This enables us to apply the Lions-Aubin lemma to derive the global convergence of $\rho_\varepsilon$ to $\rho$ on $\mathbb{R}^d$. This argument, which yields a global (not just local) convergence result, is adapted from Sugiyama and Kunii \cite{SK}. However, it is clear that when the potential under consideration shifts from the Newtonian type to the Riesz type, the condition $0<\gamma \leq d-2$ is no longer sufficient for our purposes as $d\geq 5$, thus leading to the range specified for \eqref{gamma}. Furthermore, unlike the method used in \cite{BD}, our analysis does not require the initial data to have $L^\infty$ boundedness. However, in case one considers $L^\infty$ solutions, the blow-up in $L^m$-norm \eqref{r55} implies the blow-up in $L^\infty$-norm due to the mass conservation.
	\end{rem}

	Notice that the assumptions \eqref{c1} and \eqref{c2} in the above two results have totally different formulation, we can verify that they are equivalent. These two assumptions have also been mentioned in the case $\gamma=d-2$ in \cite{CW} and \cite{KNO} separately. However, to the authors' knowledge, there is no discussion on the relation between them, even for $\gamma=d-2$.
	In the following statement, we proceed with detailed computations to show that they are in fact equivalent.
	\begin{theorem} \label{th3}
		The conditions \eqref{c1} and \eqref{c2} are equivalent. Furthermore, two sharp criteria in Theorem \eqref{th1} and Theorem \eqref{th2} are equivalent.
	\end{theorem}

	This paper is structured as follows: In Section 2, we establish the global existence of solutions stated in Theorem \ref{th1} (1) through the $L^{m_*}$-norm of the initial data, and blow-up in finite time of solutions in Theorem \ref{th1} (2) will be given. In Section 3, we prove the global existence of solutions stated in Theorem \ref{th2} through the $L^m$-norm of the initial data. We notice that the blow-up proof of Theorem \ref{th2} is omitted, since it can be similarly obtained as in \cite{BD}. In Section 4, we focus on proving the equivalence of two sharp criteria.
	
	\section{Proof of Theorem \ref{th1}.}
	This section aims to give the proof of the first main result.
	\subsection{Existence of weak solutions.}
	In this subsection, we establish the existence of weak solutions to \eqref{a1}. We follow the approach of \cite{BCL,CW,S,SK} and consider the following regularized system:
	\begin{equation}\label{r1}
		\begin{cases}
			\partial_t \rho_{\varepsilon}=\Delta[(\rho_\varepsilon+\varepsilon)^m-\varepsilon^m]-\text{div}((\rho_\varepsilon+\varepsilon)\nabla c_{\varepsilon}), &x\in \mathbb{R}^d, t\geq 0,\\
			c_{\varepsilon}=\bar{c}_d \frac{1}{(|x|^2+\varepsilon^2)^{\frac{\gamma}{2}}}\ast \rho_\varepsilon,  &x\in \mathbb{R}^d, t\geq 0,\\
			\rho_\varepsilon(x,0)=\rho_{0\varepsilon}(x), &x\in \mathbb{R}^d,
		\end{cases}
	\end{equation}
	for $\varepsilon>0$.  The non-negative initial data $\rho_{0\varepsilon}$ is the regularization of the function $\rho_0$ and it satisfies
	$$
 \|\rho_{0\varepsilon}\|_{L^1(\mathbb{R}^d)}=\|\rho_{0}\|_{L^1(\mathbb{R}^d)}.
	$$
	Moreover, we have
	\begin{equation}\label{na}
		\int_{\mathbb{R}^d} |x|^2 \rho_{0\varepsilon} \, dx \to \int_{\mathbb{R}^d} |x|^2 \rho_0 dx \quad \text{and} \quad \int_{\mathbb{R}^d}|\nabla \rho_{0\varepsilon}^m|^2dx
		\to \int_{\mathbb{R}^d}|\nabla \rho_0^m|^2dx, \quad \varepsilon\to 0.
	\end{equation}
	If $\rho_0 \in L^p$ for some $p > 1$, then $\|\rho_{0\varepsilon} - \rho_0\|_{L^p(\mathbb{R}^d)} \to 0$ as $\varepsilon \to 0$.
	
	By classical parabolic theory, the above regularized problem admits a global smooth non-negative solution $\rho_\varepsilon$ for $t>0$,  which conserves the total mass. This solution also possesses an energy functional defined by
	\begin{equation}\label{r2}
		\begin{aligned}
			\mathcal{F}_\varepsilon(\rho_\varepsilon) =& \frac{1}{m-1} \int_{\mathbb{R}^d} \bigl((\rho_\varepsilon + \varepsilon)^m - \varepsilon^m\bigr) dx - \frac{1}{2} \int_{\mathbb{R}^d} \rho_\varepsilon c_\varepsilon dx\\
			=& \frac{1}{m-1} \int_{\mathbb{R}^d} \bigl((\rho_\varepsilon + \varepsilon)^m - \varepsilon^m\bigr) dx
			-\frac{\bar{c}_d}{2} \iint_{\mathbb{R}^d \times \mathbb{R}^d} \frac{\rho_\varepsilon(x,t) \rho_\varepsilon(y,t)}{(|x-y|^2 + \varepsilon^2)^{\frac{\gamma}{2}}} dx dy.
		\end{aligned}
	\end{equation}
	The energy functional satisfies the following  entropy-entropy production relation
	\begin{equation}\label{r3}
		\frac{d}{dt} \mathcal{F}_\varepsilon(\rho_\varepsilon(\cdot, t)) + \int_{\mathbb{R}^d} (\rho_\varepsilon + \varepsilon) \left| \nabla \left( \frac{m}{m-1} (\rho_\varepsilon + \varepsilon)^{m-1} - c_\varepsilon \right) \right|^2 dx = 0.
	\end{equation}
	This leads to the monotonic decrease of the free energy in time, along with the non-negativity of entropy production. In addition, the initial free energy satisfies that as $\varepsilon\to 0$
	$$
	\mathcal{F}_\varepsilon (\rho_{0\varepsilon}) \to \mathcal{F}(\rho_0)=\frac{1}{m-1} \int_{\mathbb{R}^d} \rho_0^m(x,t) dx - \frac{1}{2} \int_{\mathbb{R}^d} \rho_0(x,t) c_0(x,t) dx,
	$$
	with $c_0=\bar{c}_d\frac{1}{|x|^{\gamma}} \ast \rho_0$.
	
	The free energy can be divided into two parts by inserting terms in the HLS inequality \eqref{hls2}
	$$
	\begin{aligned}
		\mathcal{F}_\varepsilon (\rho_\varepsilon)
		=& \frac{1}{m-1} \int_{\mathbb{R}^d} \left( (\rho_\varepsilon + \varepsilon)^m - \varepsilon^m \right) dx
		- \frac{C_{\text{HLS}}\bar{c}_d}{2} \|\rho_\varepsilon\|_{L^{m_*}(\mathbb{R}^d)}^2 \\
		&+\frac{C_{\text{HLS}}\bar{c}_d}{2} \|\rho_\varepsilon\|_{L^{m_*}(\mathbb{R}^d)}^2
		- \frac{\bar{c}_d}{2} \iint_{\mathbb{R}^d \times \mathbb{R}^d} \frac{\rho_\varepsilon(x,t) \rho_\varepsilon(y,t)}{(|x-y|^2 + \varepsilon^2)^{\frac{\gamma}{2}}} dx dy \\
		\geq & \left(\frac{1}{m-1} \int_{\mathbb{R}^d} \rho_\varepsilon^m dx
		- \frac{C_{\text{HLS}}\bar{c}_d}{2} \|\rho_\varepsilon\|_{L^{m_*}(\mathbb{R}^d)}^2\right) \\
		&+ \left(\frac{C_{\text{HLS}}\bar{c}_d}{2} \|\rho_\varepsilon\|_{L^{m_*}(\mathbb{R}^d)}^2
		- \frac{\bar{c}_d}{2} \iint_{\mathbb{R}^d \times \mathbb{R}^d} \frac{\rho_\varepsilon(x,t) \rho_\varepsilon(y,t)}{|x-y|^\gamma} dxdy\right) \\
		=:& \mathcal{F}_1 (\rho_\varepsilon) + \mathcal{F}_2 (\rho_\varepsilon).
	\end{aligned}
	$$
	We obtain from \eqref{hls2} that $\mathcal{F}_2 (\rho_\varepsilon)\geq 0$.
	The first part of the free energy from \eqref{int} is
	$$
	\begin{aligned}
		\mathcal{F}_1(\rho_\varepsilon) =& \frac{1}{m-1} \int_{\mathbb{R}^d} \rho_\varepsilon^m(x,t) \, dx - \frac{C_{\text{HLS}}\bar{c}_d}{2} \|\rho_\varepsilon\|_{L^{m_*}(\mathbb{R}^d)}^2 \\
		\geq &  \frac{1}{m-1} \|\rho_\varepsilon\|_{L^1(\mathbb{R}^d)}^{\frac{2d-m(2d-\gamma)}{\gamma}} \|\rho_\varepsilon\|_{L^{m_*}(\mathbb{R}^d)}^{\frac{2d(m-1)}{\gamma}}
		-\frac{C_{\text{HLS}}\bar{c}_d}{2} \|\rho_\varepsilon\|_{L^{m_*}(\mathbb{R}^d)}^2.
	\end{aligned}
	$$
	Let
	$$
	f_1(s) = \frac{1}{m-1} \|\rho_0\|_{L^1(\mathbb{R}^d)}^{\frac{2d-m(2d-\gamma)}{\gamma}} s - \frac{C_{\text{HLS}}\bar{c}_d}{2} s^{\frac{\gamma}{d(m-1)}}.
	$$
	Thus, we have a lower bound of the first part of free energy, i.e.,
	$$
	f_1\left( \|\rho_\varepsilon\|_{L^{m_*}(\mathbb{R}^d)}^{\frac{2d(m-1)}{\gamma}} \right) \leq \mathcal{F}_1(\rho_\varepsilon).
	$$
	
	\begin{lemma}\label{le3}
		If the initial free energy satisfies \eqref{c1}
		$$
		\mathcal{F}(\rho_{0\varepsilon})<\frac{\gamma-d(m-1)}{\gamma(m-1)} \left(\frac{2d^2\alpha(d)}{ C_{\text{HLS}}} \right)^{\frac{d(m-1)}{\gamma-d(m-1)}} \|\rho_{0\varepsilon}\|_{L^1(\mathbb{R}^d)}^{\frac{2d-m(2d-\gamma)}{\gamma-d(m-1)}}
		$$
		and $\|\rho_{0\varepsilon}\|_{L^{m_*}(\mathbb{R}^d)} < s_*$ with $m\in (m_*,m^*)$ and $0<\gamma \leq d-2$, and
		$\rho_\varepsilon$ is a solution of problem \eqref{r1}, then there exists a constant $\mu_1 < 1$ such that
		\begin{equation}\label{r4}
			\sup_{t\in (0,\infty)}\|\rho_\varepsilon(\cdot, t)\|_{L^{m_*}(\mathbb{R}^d)} < \mu_1 s_*,
		\end{equation}
		where $s_*$ is given in \eqref{s}.
	\end{lemma}
	
	\begin{proof}
		The fact $\frac{2d}{2d-\gamma}< m < \frac{d+\gamma}{d}$ leads to  $\frac{\gamma}{d(m-1)} > 1$. Thus, we obtain that $f_1(s)$ is strictly concave for $0 < s < \infty$. Direct calculation shows that
		$$
		f_1'(s) = \frac{1}{m-1} \|\rho_0\|_{L^1(\mathbb{R}^d)}^{\frac{2d-m(2d-\gamma)}{\gamma}}
		-\frac{C_{\text{HLS}}\bar{c}_d}{2} s^{\frac{\gamma}{d(m-1)}-1},
		$$
		which implies that
		$
		\tilde{s}_1 = \left( \frac{2d^2 \alpha(d) }{C_{\text{HLS}}} \|\rho_0\|_{L^1(\mathbb{R}^d)}^{\frac{\gamma m-2d(m-1)}{\gamma}}\right)^{\frac{d(m-1)}{\gamma-d(m-1)}}
		$
		is the unique maximum point of  $f_1(s)$. Therefore, $f_1(s)$ is monotone increasing for $0 < s < \tilde{s}_1$.
		In the case of $ \mathcal{F}_\varepsilon (\rho_{0\varepsilon}) < f_1(\tilde{s}_1) $, we can in fact choose $0< \delta < 1$ such that $ \mathcal{F}_\varepsilon (\rho_{0\varepsilon}) < \delta f_1(\tilde{s}_1) $.
		Combining this with the initial conditions, \eqref{hls2} and the monotonicity of the free energy, we have
		$$
		f_1\left( \|\rho_\varepsilon\|_{L^{m_*}(\mathbb{R}^d)}^{\frac{2d(m-1)}{\gamma}} \right) \leq \mathcal{F}_1(\rho_\varepsilon) \leq \mathcal{F}_\varepsilon (\rho_\varepsilon) \leq \mathcal{F}_\varepsilon (\rho_{0\varepsilon}) < \delta f_1(\tilde{s}_1).
		$$
		If the initial data satisfies  $\|\rho_{0\varepsilon}\|_{L^{m_*}(\mathbb{R}^d)}^{\frac{2d(m-1)}{\gamma}}< \tilde{s}_1$, the fact that $f_1(s)$ is increasing on $(0, \tilde{s}_1)$ implies there exists $\mu_1 < 1$ such that
		$$
		\|\rho_\varepsilon\|_{L^{m_*}(\mathbb{R}^d)}
		<\mu_1\tilde{s}_1^{\frac{\gamma}{2d(m-1)}}, \quad \text{for all } t> 0,
		$$
		which implies \eqref{r4}.
	\end{proof}
	
	\begin{lemma}\label{le4}
		Assume $\rho_{0 \varepsilon} \in L_+^{1}(\mathbb{R}^d) \cap L^{p}(\mathbb{R}^d)$ with $p>1$, $\|\rho_{0\varepsilon}\|_{L^{m_*}(\mathbb{R}^d)} < s_*$ with $m\in (m_*, m^*)$ and $0<\gamma \leq d-2$,  and  the initial free energy satisfies \eqref{c1}. Let $\rho_{\varepsilon}$ be a smooth solution of the regularized system \eqref{r1}, then
		\begin{equation}\label{r6}
			\|\rho_{\varepsilon}\|_{L^{\infty}(0,\infty;L^{p}(\mathbb{R}^d))} \leq C\left(d, \gamma, p, m, \mu_1, \tilde{s}_1, \|\rho_0\|_{L^1(\mathbb{R}^d)}, \|\rho_{0\varepsilon}\|_{L^p(\mathbb{R}^d)}\right),
		\end{equation}
		and for $T>0$, $r_1=\frac{d(p+1)}{2d-\gamma-2}$, it holds
		\begin{equation}\label{r61}
			\begin{aligned}
				&\|\rho_{\varepsilon}\|_{ L^{p+1}(0,T;L^{r_1}(\mathbb{R}^d))}+	\|\nabla\rho_{\varepsilon}^{\frac{m+p-1}{2}}\|_{L^{2}(0,T;L^{2}(\mathbb{R}^d))}\\
				 \leq & C\left(d, \gamma, p, m, \mu_1, \tilde{s}_1,  \|\rho_0\|_{L^1(\mathbb{R}^d)}, \|{\rho_{0\varepsilon}}\|_{L^p(\mathbb{R}^d)}, T\right).
			\end{aligned}
		\end{equation}
		Furthermore, for $\alpha \in \left(\frac{d}{\gamma+1}, \frac{pd}{p(\gamma+1-d)+d}\right]$, it holds
		\begin{equation}\label{r7}
			\|\nabla c_{\varepsilon}\|_{L^{\infty}(0,\infty;L^{\alpha}(\mathbb{R}^d))} \leq  C\left(d,\gamma, p, \mu_1, \tilde{s}_1,  \|\rho_0\|_{L^1(\mathbb{R}^d)},  \|\rho_{0\varepsilon}\|_{L^p(\mathbb{R}^d)}\right).
		\end{equation}
	\end{lemma}
	
	\begin{proof}
		For convenience, we use $L^p$ to represent $L^p(\mathbb{R}^d)$. Multiplying the first equation of \eqref{r1} by $p \rho_{\varepsilon}^{p-1}$ with $p > 1$, we have
		\begin{equation}\label{r11}
			\begin{aligned}
				\frac{d}{dt} \int_{\mathbb{R}^d} \rho_{\varepsilon}^{p} dx
				\leq & -pm(p-1) \int_{\mathbb{R}^d} \rho_{\varepsilon}^{p+m-3} |\nabla \rho_{\varepsilon}|^{2} dx
				+ (p-1) \int_{\mathbb{R}^d} \rho_{\varepsilon}^p(-\Delta c_{\varepsilon}) dx \\
				&+ p \varepsilon \int_{\mathbb{R}^d} \rho_{\varepsilon}^{p-1}(-\Delta c_{\varepsilon}) dx
				=:H_1+H_2+H_3.\\
			\end{aligned}
		\end{equation}
		For $H_1$, through simple calculations, we obtain
		\begin{equation}\label{r12}
			H_1=-\frac{4pm(p-1)}{(m+p-1)^{2}} \int_{\mathbb{R}^d} |\nabla \rho_{\varepsilon}^{\frac{m+p-1}{2}}|^{2} dx.
		\end{equation}
		According to
		\begin{equation}\label{r28}
			-\Delta c_{\varepsilon}
			\leq \bar{c}_dd\gamma\frac{1}{(|x|^2+\varepsilon^2)^{\frac{\gamma+2}{2}}}\ast\rho_{\varepsilon}
			\leq \bar{c}_dd\gamma\frac{1}{|x|^{\gamma+2}}\ast\rho_{\varepsilon},
		\end{equation}
		we obtain from \eqref{hls} that
		$$
		H_2\leq (p-1)\bar{c}_dd\gamma \int_{\mathbb{R}^d}\int_{\mathbb{R}^d}\frac{\rho_{\varepsilon}^p(x)\rho_{\varepsilon}(y)}{|x-y|^{\gamma+2}}dxdy
		\leq  (p-1)\bar{c}_dd\gamma C_{\text{HLS}}\|\rho_{\varepsilon}\|_{L^{r_1}}^{p+1}.
		$$
		Here $\frac{p}{r_1}+\frac{1}{r_1}+\frac{\gamma+2}{d}=2$, i.e., $r_1=\frac{d(p+1)}{2d-\gamma-2}$.
		According to \eqref{GNS}, we have
		\begin{equation}\label{r13}
			\|\rho_{\varepsilon}\|_{L^{r_1}}^{p+1}
			=\|\rho_{\varepsilon}^{\frac{m+p-1}{2}}\|_{L^{\frac{2r_1}{m+p-1}}}^{\frac{2(p+1)}{m+p-1}}
			\leq S_d^{-\frac{(p+1)\theta_1}{m+p-1}}\|\nabla \rho_{\varepsilon}^{\frac{m+p-1}{2}}\|_{L^2}^{\frac{2(p+1)\theta_1}{m+p-1}}\|\rho_{\varepsilon}\|_{L^{m_*}}^{(1-\theta_1)(p+1)}.
		\end{equation}
		Here $\theta_1=\frac{(2d-\gamma)(m+p-1)r_1-2d(m+p-1)}{(2d-\gamma)(m+p-1)r_1-2r_1(d-2)}$ with $\frac{2d}{2d-\gamma}<r_1<\frac{d(m+p-1)}{d-2}$.
		Then, we have
		$$
		\frac{2(p+1)\theta_1}{m+p-1}=2+\frac{2(2d-\gamma)(2-m)-4(d-\gamma)}{(2d-\gamma)(m+p-1)-2(d-2)}<2.
		$$
		The fact $m>m_*$ implies $\frac{(2d-\gamma)(2-m)}{2d}-\frac{d-\gamma}{d}<0$. Thus, we use Young's inequality to obtain
		\begin{equation}\label{r9}
			H_2\leq \frac{pm(p-1)}{(m+p-1)^2}\|\nabla \rho_{\varepsilon}^{\frac{m+p-1}{2}}\|_{L^2}^2
			+C\left(\gamma, p, m, d\right) \|\rho_{\varepsilon}\|_{L^{m_*}}^{(1-\theta_1)(p+1)p_1},
		\end{equation}		
		where $p_1=1+\frac{(p+1)\theta_1}{m+p-1-(p+1)\theta_1}$.
		By the same method, we obtain
		\begin{equation}\label{r10}
			H_3\leq \frac{pm(p-1)}{(m+p-1)^2} \|\nabla \rho_{\varepsilon}^{\frac{m+p-1}{2}}\|_{L^2}^2
			+C\left(\gamma, p, m, d \right) \|\rho_{\varepsilon}\|_{L^{m_*}}^{p_2},
		\end{equation}
		where $p_2=\frac{2(m+p-1)(2d-\gamma-2)-2p(d-2)}{(m+1)(2d-\gamma)-2d}>1$.
		Inserting \eqref{r12}, \eqref{r9} and \eqref{r10} into \eqref{r11}, we deduce from \eqref{r4} that
		\begin{equation}\label{r56}
			\begin{aligned}
				\frac{d}{dt} \int_{\mathbb{R}^d} \rho_{\varepsilon}^{p} dx
				\leq & -\frac{2pm(p-1)}{(m+p-1)^2}\|\nabla \rho_{\varepsilon}^{\frac{m+p-1}{2}}\|_{L^2}^2\\
				&+ C\left(d, \gamma, p, m\right)\left( \|\rho_\varepsilon\|_{L^{m_*}}^{(1-\theta_1)(p+1)p_1}
				+\|\rho_\varepsilon\|_{L^{m_*}}^{p_2}\right).
			\end{aligned}
		\end{equation}
		Applying interpolation inequality with $1<p<\frac{(m+p-1)d}{d-2}$, \eqref{So} and Young's inequality, we obtain
		$$
		\begin{aligned}
			\|\rho_\varepsilon\|_{L^{p}}^{p}
			\leq & \|\rho_\varepsilon\|_{L^1}^{p\theta_2}\left(S_d^{-\frac{1}{2}}\|\nabla \rho_\varepsilon^{\frac{m+p-1}{2}}\|_{L^2}\right)^{\frac{2p(1-\theta_2)}{m+p-1}}\\
			\leq &\frac{2pm(p-1)}{(m+p-1)^2}\|\nabla \rho_\varepsilon^{\frac{m+p-1}{2}}\|_{L^2}^{2}
			+C\left(d, \gamma, p, m \right)\|\rho_\varepsilon\|_{L^1}^{p_3}.
		\end{aligned}
		$$
		Here $\theta_2=\frac{d(m-1)+2p}{p[d(m+p-1)+d-2]}$, which leads to $\frac{2p(1-\theta_2)}{m+p-1}<2$, and $p_3=\frac{(d(m-1)-2p)(m+p-1)}{d(m+p-1)(m-1)-2m+2p+2}$.
		Thus, combining above two results, we can derive
		$$
		\frac{d}{dt} \int_{\mathbb{R}^d} \rho_{\varepsilon}^{p} dx
		\leq -\|\rho_\varepsilon\|_{L^{p}}^{p}+C\left(d, \gamma, p, m\right)\left( \|\rho_\varepsilon\|_{L^{m_*}}^{(1-\theta_1)(p+1)p_1}
		+\|\rho_\varepsilon\|_{L^{m_*}}^{p_2}+\|\rho_\varepsilon\|_{L^1}^{p_3}\right).
		$$
		Thus, we have from \eqref{r4} and $\rho_\varepsilon\in L^\infty(0,\infty; L^1)$ that, for all $t>0$,
		$$
		\|\rho_\varepsilon\|_{L^p}^p\leq \|\rho_{0\varepsilon}\|_{L^p}^p+C\left(d, \gamma, p, m, \mu_1, \tilde{s}_1,  \|\rho_0\|_{L^1}\right),
		$$
		i.e.,
		
		$$
		\|\rho_{\varepsilon}\|_{L^{\infty}(0,\infty;L^p)}\leq C\left(d, \gamma, p, m, \mu_1, \tilde{s}_1, \|\rho_0\|_{L^1},  \|\rho_{0\varepsilon}\|_{L^p}\right).
		$$
		Taking integral with respect to the time variable of \eqref{r56}, we obtain $\forall T>0$ that
		$$
		\begin{aligned}
			&\sup_{0\leq t\leq  T}\int_{\mathbb{R}^d}\rho_{\varepsilon}^p(x,t)dx
			+\frac{2pm(p-1)}{(m+p-1)^2}\int_0^T\int_{\mathbb{R}^d}|\nabla \rho_{\varepsilon}^{\frac{m+p-1}{2}}|^2dxdt\\
			\leq & C\left(d, \gamma, p, m, \mu_1, \tilde{s}_1, \|\rho_0\|_{L^1}, \|{\rho_{0\varepsilon}}\|_{L^p}, T\right).
		\end{aligned}
		$$
		Thus, we deduce that
		\begin{equation}\label{r33}
			\|\nabla\rho_{\varepsilon}^{\frac{m+p-1}{2}}\|_{L^{2}(0,T;L^{2})} \leq  C\left(d, \gamma, p, m, \mu_1, \tilde{s}_1, \|\rho_0\|_{L^1}, \|{\rho_{0\varepsilon}}\|_{L^p}, T\right).
		\end{equation}
		Combining this result with \eqref{r13} and \eqref{r33}, we deduce that
		$$
		\|\rho_\varepsilon\|_{L^{p+1}(0,T;L^{r_1})}\leq  C\left(d, \gamma, p, m, \mu_1, \tilde{s}_1, \|\rho_0\|_{L^1}, \|{\rho_{0\varepsilon}}\|_{L^p}, T\right),
		\quad r_1=\frac{d(p+1)}{2d-\gamma-2}.
		$$
		Finally, the estimates for $\nabla c_\varepsilon$ can be obtained by applying the weak Young inequality (see \cite[Page 107]{LL}) to
		$$
		\nabla c_{\varepsilon}=-\gamma \bar{c}_d\frac{x}{(|x|^2+\varepsilon^2)^{\frac{\gamma+2}{2}}}\ast \rho_\varepsilon.
		$$
		Namely,
		$$
		\begin{aligned}
			\| \nabla c_{\varepsilon} \|_{L^\alpha} = &\gamma \bar{c}_d \left\| \frac{x}{(|x|^2+\varepsilon^2)^{\frac{\gamma+2}{2}}} * \rho_\varepsilon \right\|_{L^\alpha}
			\leq \gamma \bar{c}_d \left\| \frac{x}{|x|^{\gamma+2}} \right\|_{L^q_w} \| \rho_{\varepsilon} \|_{L^{q'}} \\
			\leq &C \| \rho_{\varepsilon} \|_{L^{q'}} \leq C\left(d,\gamma, p, \mu_1, \tilde{s}_1, \|\rho_0\|_{L^1}, \|\rho_{0\varepsilon}\|_{L^p}\right),
		\end{aligned}
		$$
		where $q' \in (1, p]$, $q = \frac{d}{\gamma+1}$ and $\frac{1}{\alpha} + 1 = \frac{1}{q} + \frac{1}{q'}$. So we have $\alpha \in \left( \frac{d}{\gamma+1}, \frac{pd}{p(\gamma+1-d)+d} \right]$.
	\end{proof}

	\begin{lemma}\label{le5}
		Assume the assumptions of Lemma \ref{le4} hold with $p=m$, then it holds that for any $T>0$
		\begin{equation}\label{r14}
			\|\nabla \rho_{\varepsilon}\|_{L^2(0,T; L^{\frac{2m}{3-m}}(\mathbb{R}^d))} \leq
			C\left(d, \gamma, m, \mu_1, \tilde{s}_1,  \|\rho_0\|_{L^1(\mathbb{R}^d)},  \|{\rho_{0\varepsilon}}\|_{L^m(\mathbb{R}^d)}, T\right), \quad \text{for } m < \frac{3}{2},
		\end{equation}
		
		\begin{equation}\label{r15}
			\| \nabla \rho_{\varepsilon} \|_{L^2(0,T; L^{2}(\mathbb{R}^d))} \leq 	C\left(d, \gamma, m, \mu_1, \tilde{s}_1,  \|\rho_0\|_{L^1(\mathbb{R}^d)},  \|{\rho_{0\varepsilon}}\|_{L^m(\mathbb{R}^d)}, T\right), \quad \text{for } m \geq \frac{3}{2}.
		\end{equation}
		
	\end{lemma}
	
	\begin{proof}
		For the case of $m <\frac{3}{2}$, we refer to the proof of \eqref{r14} to \cite[Lemma 5.5]{WCH}.

		For the case of $m \geq \frac{3}{2}$, taking $\rho_{\varepsilon}^{2-m}$ as test function in \eqref{r1}, we obtain
		\begin{equation}\label{r22}
			\begin{aligned}
				&\frac{1}{3-m} \frac{d}{dt} \int_{\mathbb{R}^d} \rho_{\varepsilon}^{3-m}dx + m(2-m) \int_{\mathbb{R}^d} |\nabla \rho_{\varepsilon}|^2 dx\\
				\leq & -\frac{2-m}{3-m} \int_{\mathbb{R}^d} \rho_{\varepsilon}^{3-m} \Delta c_{\varepsilon} dx -\varepsilon \int_{\mathbb{R}^d} \rho_{\varepsilon}^{2-m} \Delta c_{\varepsilon}dx:=H_4+H_5.
			\end{aligned}
		\end{equation}
		Clearly, it follows from \eqref{hls} that
		\begin{equation}\label{r18}
			H_4
			\leq\frac{2-m}{3-m}\bar{c}_dd\gamma \|\rho_{\varepsilon}\|_{L^{r_3}}^{4-m},
		\end{equation}
		where $\frac{3-m}{r_3}+\frac{1}{r_3}+\frac{\gamma+2}{d}=2$ leads to $r_3=\frac{(4-m)d}{2d-\gamma-2}$. Then, we apply the interpolation inequality because of $1<r_3\leq r_1$ from $m\geq \frac{3}{2}$, and obtain
		\begin{equation}\label{r19}
			\|\rho_{\varepsilon}\|_{L^{r_3}}^{4-m}
			\leq \|\rho_{\varepsilon}\|_{L^1}^{\theta_3(4-m)}
			\|\rho_{\varepsilon}\|_{L^{r_1}}^{(1-\theta_3)(4-m)}, \quad \theta_3=\frac{(2d-\gamma-2)(2m-3)}{(4-m)[d(m+1)-(2d-\gamma-2)]}.
		\end{equation}
		Next, there holds from \eqref{GNS} that
		\begin{equation}\label{r20}
			\begin{aligned}
				\|\rho_{\varepsilon}\|_{L^{r_1}}^{(1-\theta_3)(4-m)}
				=&\|\rho_{\varepsilon}^{m-\frac{1}{2}}\|_{L^{\frac{r_1}{m-\frac{1}{2}}}}^{\frac{(1-\theta_3)(4-m)}{m-\frac{1}{2}}}\\
				\leq & S_d^{-\frac{(1-\theta_3)(4-m)}{m+1}}
				\|\nabla \rho_\varepsilon^{m-\frac{1}{2}}\|_{L^2}^{\frac{2(1-\theta_3)(4-m)}{m+1}}
				\|\rho_\varepsilon\|_{L^{\frac{d(2-m)}{d-\gamma}}}^{(1-\theta_3)(4-m)\frac{2-m}{m+1}}.
			\end{aligned}
		\end{equation}
		By the same argument, we also have
		\begin{equation}\label{r23}
			H_5\leq \varepsilon \bar{c}_dd\gamma S_d^{-\frac{(1-\theta_4)(3-m)}{m+1}}\|\rho_\varepsilon\|_{L^1}^{\theta_4(3-m)}
			\|\nabla \rho_{\varepsilon}^{m-\frac{1}{2}}\|_{L^2}^{\frac{2(1-\theta_4)(3-m)}{m+1}}
			\|\rho_\varepsilon\|_{L^{\frac{d(2-m)}{d-\gamma}}}^{(1-\theta_4)(3-m)\frac{(2-m)}{m+1}}.
		\end{equation}
		Here $\theta_4=\frac{(2d-\gamma-2)(2m-2)}{(3-m)[d(m+1)-(2d-\gamma-2)]}$.
		Inserting \eqref{r18}-\eqref{r23} into \eqref{r22}, taking integral with respect to the time variable and using Young's inequality and the interpolation inequality with $1<3-m< m$ (if $3-m=m$, the first term on the right directly depends on $\|\rho_{0\varepsilon}\|_{L^m}$ and does not require interpolation inequality), we obtain that for any $T>0$
		$$
		\begin{aligned}
			&\frac{1}{3-m} \int_{\mathbb{R}^d} \rho_{\varepsilon}^{3-m}dx + m(2-m)\int_0^T \int_{\mathbb{R}^d} |\nabla \rho_{\varepsilon}|^2 dxdt\\
			\leq & \frac{1}{3-m}\|\rho_{0\varepsilon}\|_{L^{3-m}}^{3-m}+	C\left(d, \gamma,m, \|\rho_\varepsilon\|_{L^1}, \mu_1, \tilde{s}_1, T, \|{\rho_{0\varepsilon}}\|_{L^m}\right)\\
			\leq & 	C\left(d, \gamma, m, \mu_1, \tilde{s}_1, \|\rho_0\|_{L^1}, \|{\rho_{0\varepsilon}}\|_{L^m}, T\right).
		\end{aligned}
		$$
		Then the fact $\frac{2(1-\theta_4)(3-m)}{m+1}<\frac{2(1-\theta_3)(4-m)}{m+1}<2$ and $1<\frac{d(2-m)}{d-\gamma}<m$ means that $\eqref{r15}$ holds.
	\end{proof}
	
	\begin{lemma}\label{le6}
		Assume the assumptions of Lemma \ref{le3} hold with $p=m$, then it holds that for any $T>0$,
		$$
		\partial_t \rho_\varepsilon \in L^2\left(0,T;W_{loc}^{-1,s}(\mathbb{R}^d)\right), \quad s\in \min\left\{\frac{2m}{m+1}, \frac{m(m+1)d}{(2d-\gamma-2)m+[m(\gamma+1-d)+d](m+1)}\right\}.
		$$
	\end{lemma}
	\begin{proof}
		By using the weak formulation of the equation, we know the estimate for time derivative $\partial_t \rho_\varepsilon$ can be obtained directly from the estimates on $\nabla[(\rho_\varepsilon+\varepsilon)^m-\varepsilon^m]$ and $(\rho_\varepsilon+\varepsilon)\nabla c_\varepsilon$.
		
		In fact, there holds
		$$
		|\nabla[(\rho_\varepsilon+\varepsilon)^m-\varepsilon^m]|
		=m(\rho_{\varepsilon}+\varepsilon)^{m-1}\nabla \rho_\varepsilon
		\leq  \nabla \rho_\varepsilon^m+m\varepsilon^{m-1}\nabla \rho_\varepsilon.
		$$
		The relationship $\nabla \rho_\varepsilon^m=\frac{2m}{2m-1}\rho_\varepsilon^{\frac{1}{2}}\nabla \rho_\varepsilon^{m-\frac{1}{2}}$ and the H\"{o}lder inequality yield
		$$
		\int_{\mathbb{R}^d}|\nabla \rho_\varepsilon^m|^{\frac{2m}{m+1}}dx
		\leq \left(\int_{\mathbb{R}^d}\rho_\varepsilon^mdx\right)^{\frac{1}{m+1}}
		\left(\int_{\mathbb{R}^d}|\nabla \rho_\varepsilon^{m-\frac{1}{2}}|^2dx\right)^{\frac{m}{m+1}},
		$$
		which implies from \eqref{r6} that for any $T>0$
		$$
		\int_0^T\|\nabla \rho_\varepsilon^m\|_{L^{\frac{2m}{m+1}}}^2dt
		\leq \int_0^T\|\rho_\varepsilon\|_{L^m}\|\nabla \rho_\varepsilon^{m-\frac{1}{2}}\|_{L^2}^2dt
		\leq C\left(d, \gamma, m, \mu_1, \tilde{s}_1, \|\rho_0\|_{L^1},  \|{\rho_{0\varepsilon}}\|_{L^m}, T\right).
		$$
		Namely,
		$$
		\|\nabla \rho_\varepsilon^m\|_ {L^2\left(0,T;L^{\frac{2m}{m+1}}\right)}
		\leq C\left(d, \gamma, m, \mu_1, \tilde{s}_1, \|\rho_0\|_{L^1}, T, \|{\rho_{0\varepsilon}}\|_{L^m}\right).
		$$
		Combining with Lemma \ref{le5} and $\frac{2m}{m+1}<\min\{2,\frac{2m}{3-m}\}$, we obtain that
		\begin{equation}\label{r25}
			 \|\nabla[(\rho_\varepsilon+\varepsilon)^m-\varepsilon^m]
			 \|_{L^2\left(0,T;L_{loc}^{\frac{2m}{m+1}}(\mathbb{R}^d)\right)}\leq C\left(d, \gamma, m, \mu_1, \tilde{s}_1, \|\rho_0\|_{L^1}, \|{\rho_{0\varepsilon}}\|_{L^m}, T\right).
		\end{equation}
		Applying the H\"{o}lder inequality, we have
		$$
		\|\rho_\varepsilon\nabla c_\varepsilon\|_{L^{\frac{m(m+1)d}{(2d-\gamma-2)m+[m(\gamma+1-d)+d](m+1)}}}
		\leq \|\rho_\varepsilon\|_{L^{r_1}}\|\nabla c_{\varepsilon}\|_{L^{\frac{md}{m(\gamma+1-d)}}}.
		$$
		According to Lemma \ref{le4} and the result above with
		$$
		\frac{m(m+1)d}{(2d-\gamma-2)m+[m(\gamma+1-d)+d](m+1)}\in \left(\frac{d}{\gamma+1}, \frac{md}{m(\gamma+1-d)+d}\right],
		$$
		we have
		\begin{equation}\label{r26}
			\begin{aligned}
				&\|(\rho_\varepsilon+\varepsilon)\nabla c_\varepsilon\|_{L^{m+1}\left(0,T;
					L^\frac{m(m+1)d}{(2d-\gamma-2)m+[m(\gamma+1-d)+d](m+1)}\right)}\\
				\leq & C\left(d, \gamma, m, \mu_1, \tilde{s}_1, \|\rho_0\|_{L^1}, \|{\rho_{0\varepsilon}}\|_{L^m}, T\right).
			\end{aligned}
		\end{equation}
		Obviously, \eqref{r25} and \eqref{r26} lead to
		\begin{equation}\label{r27}
			\partial_t \rho_\varepsilon \in L^2\left(0,T;W_{loc}^{-1,s}(\mathbb{R}^d)\right).
		\end{equation}
		Here $s\in \min\{\frac{2m}{m+1}, \frac{m(m+1)d}{(2d-\gamma-2)m+[m(\gamma+1-d)+d](m+1)}\}>1$.
	\end{proof}
	
	\begin{proof}[Proof of Theorem \ref{th1} (1)]
		For $m<\frac{3}{2}$, we deduce from \eqref{r14} that
		\begin{equation}\label{r17}
			\rho_{\varepsilon}\in L^2(0,T;W^{1,\frac{2m}{3-m}}(\mathbb{R}^d)),\quad \text{for any fixed } T.
		\end{equation}
		For $m\geq \frac{3}{2}$,  we use \eqref{So} to obtain
		$$
		\|\rho_\varepsilon\|_{L^{\frac{2d}{d-2}}}\leq S_d^{-\frac{1}{2}}\|\nabla \rho_\varepsilon\|_{L^{2}}.
		$$
		Noticing $\frac{2d}{d-2}\geq 2$ and combining with the fact that  $\rho_{\varepsilon} \in L^{\infty}(0,\infty;L^{m}) \cap L^{m+1}(0,T;L^{r_1})$ with $r_1=\frac{d(m+1)}{2d-\gamma-2}$, we obtain
		\begin{equation}\label{r24}
			\rho_{\varepsilon}\in L^2\left(0,T;W^{1,2}(\mathbb{R}^d)\right),\quad \text{for any } T>0.
		\end{equation}
		According to \eqref{r27} \eqref{r17} and \eqref{r24}, we apply Lions-Aubin Lemma \cite{A} (see also \cite[Chapter IV, Section 4]{Lions}) to show that there exists a subsequence $\rho_\varepsilon$ without relabeling such that for any $T>0$ and any bounded domain $\Omega\subset \mathbb{R}^d$
		$$
		\rho_\varepsilon \to \rho, \quad \text{in } L^2\left(0,T;L^{\bar{p}'}(\Omega)\right),
		$$
		where $\max\{\frac{2md}{(m+1)d+2m}, \frac{m(m+1)d}{m(m+1)(\gamma+2-d)+(m+1)d+(2d-\gamma-2)m}\}\leq \bar{p}'<\min\left\{\frac{2md}{(3-m)d-2m}, \frac{2d}{d-2}\right\}$.  Let $\{B_k\}_{k=1}^{\infty}$
		be a sequence of balls centered at $0$ with radius $R_k$, $R_k\to \infty$. By a standard diagonal argument, we can find a subsequence of $\rho_\varepsilon$. Without relabeling, we have the following
		uniform strong convergence:
		$$
		\rho_\varepsilon \to \rho, \quad \text{in } L^2\left(0,T; L^{\bar{p}'}(B_k)\right) \quad \forall k.
		$$
		Similar to the argument in \cite{BL,CLW}, we take limit of \eqref{r1} and obtain that for any test function $\psi \in C_0^\infty(\mathbb{R}^d)$ in Definition \ref{de}, there exists a global weak solution $\rho$. This completes the proof of Theorem \eqref{th1} (1).
	\end{proof}

	\subsection{Blow up of solutions.}
	In this subsection, we prove the finite-time blow-up of solutions. Prior to establishing the blow-up result, we first present a key lemma demonstrating that for subcritical initial data, the quantity $\|\rho_\varepsilon\|_{L^{m_*}}$ admits a lower bound.
	
	\begin{lemma}\label{le1}
		Let $\rho$ be a solution of problem \eqref{a1}. If the initial free energy satisfies \eqref{c1} and $\|\rho_0\|_{L^{m_*}}>s_*$, then there exists a constant $\mu_2 > 1$ such that
		$$
		\|\rho(\cdot, t)\|_{L^{m_*}} > (\mu_2 \tilde{s}_1)^{\frac{\gamma}{2d(m-1)}}, \quad \text{for all } t > 0.
		$$
	\end{lemma}
	\begin{proof}
		Fix $\delta>0$ such that
		$
		\mathcal{F}(\rho_0) < \delta f_1(\tilde{s}_1).
		$
		Taking a similar computation in Lemma \ref{le3}, we find that
		if initially $\|\rho_{0\varepsilon}\|_{L^{m_*}}^{\frac{2d(m-1)}{\gamma}} > \tilde{s}_1$, the fact that $f_1(s)$ is decreasing on $s>\tilde{s}_1$ implies there exists $\mu_2 > 1$ such that
		$$
		\|\rho_\varepsilon\|_{L^{m_*}}^{\frac{2d(m-1)}{\gamma}}
		> \mu_2\tilde{s}_1, \quad \mbox{for all } t>0.
		$$
	\end{proof}
	
	Next, we turn our attention to analyzing the time evolution of the second moment. Then we can prove of blow up in Theorem \ref{th1}.
	\begin{proof}
		A direct computation shows that the time derivative of the second moment is given by
		$$
		\frac{dm_2(t)}{dt}=\frac{d}{dt}\int_{\mathbb{R}^d}|x|^2\rho(x,t)dx
		=(2d-\frac{2\gamma}{m-1})\|\rho\|_{L^m}^m+2\gamma\mathcal{F}(\rho).
		$$
		The fact $1<\frac{2d}{2d-\gamma}<m<\frac{d+\gamma}{d}$ leads to $2d-\frac{2\gamma}{m-1}<0$. Thus, using \eqref{int}, the decreasing properties of free energy and Lemma \ref{le1} with $\mu_2 > 1$, we obtain
		$$
		\begin{aligned}
			\frac{dm_2(t)}{dt}\leq &(2d-\frac{2\gamma}{m-1})\|\rho_{0\varepsilon}\|_{L^1}^{\frac{\theta-1}{\theta}m}\|\rho\|_{L^{m_*}}^{\frac{m}{\theta}}+2\gamma\mathcal{F}(\rho_0)\\
			<&(2d-\frac{2\gamma}{m-1})\|\rho_0\|_{L^1}^{\frac{\theta-1}{\theta}m}(\mu_2-1)\tilde{s}_1+\tilde{s}_1(2d\|\rho_0\|_{L^1}^{\frac{\theta-1}{\theta}m}-\gamma C_{\text{HLS}}\overline{c}_d(\tilde{s}_1)^{\frac{2\theta}{m}-1})\\
			=&(2d-\frac{2\gamma}{m-1})\|\rho_0\|_{L^1}^{\frac{\theta-1}{\theta}m}(\mu_2-1)\tilde{s}_1<0,
		\end{aligned}
		$$
		which impiles that there exists the finite time $T^*$ such that
		$
		\lim_{t\to T^*}m_2(t)=0.
		$
		Using the H\"{o}lder inequality, we have
		$$
		\int_{\mathbb{R}^d}\rho dx=\int_{B_R(0)}\rho(x)dx+\int_{\mathbb{R}^d\setminus B_R(0)}\rho(x)dx
		\leq CR^{\frac{d(m-1)}{m}}\|\rho\|_{L^m}+\frac{1}{R^2}m_2(t).
		$$
		Choosing $R=\left(\frac{m_2(t)}{C\|\rho\|_{L^m}}\right)^{\frac{m}{d(m-1)+2m}}$ such that $CR^{\frac{d(m-1)}{m}}\|\rho\|_{L^m}=\frac{1}{R^2}m_2(t)$, we deduce
		$$
		\|\rho_0\|_{L^1}
		\leq 2C\left(\frac{m_2(t)}{C\|\rho\|_{L^m}}\right)^{\frac{d(m-1)}{d(m-1)+2m}}\|\rho\|_{L^m}
		\leq 2Cm_2(t)^{\frac{d(m-1)}{d(m-1)+2m}}\|\rho\|_{L^m}^{\frac{2m}{d(m-1)+2m}}.
		$$
		Due to the fact that the mass is a fixed constant, the above inequality implies
		$$
		\limsup_{t\to T^*}\|\rho(\cdot,t)\|_{L^m}\geq
		\lim_{t\to T^*} \left(\frac{\|\rho_0\|_{L^1}}{2Cm_2(t)^{\frac{d(m-1)}{d(m-1)+2m}}}\right)^{\frac{d(m-1)+2m}{2m}}=\infty.
		$$
		Thus, the weak solution $\rho(x,t)$ blows up and this completes the proof of Theorem \ref{th1} (2).
	\end{proof}

	\section{Proof of Theorem \ref{th2}.}
	In this section, we present the proof of Theorem \ref{th2}. We omit the proof of finite-time blow-up for solutions, which can be found in \cite[Proposition 11]{BD} and \cite[Section 3.3]{CW}. We concentrate in establishing the global existence of solutions to \eqref{a1}. Several lemmas are given to deal with our problem.
	
	Under the assumptions in Theorem \ref{th2}, one can obtain a priori estimate, which shows that in the case of supercritical initial data, the quantity  $\rho_\varepsilon$ can be bounded from above. This is a crucial step in obtaining the uniform boundedness estimate
	 of $L^p$ for  $\rho_\varepsilon$. The existence result does not require $L^\infty$ boundedness for initial data, which was used in \cite{B}. Furthermore, we provide additional moment estimates to establish the compact embedding for Lions-Aubin Lemma.

	\begin{lemma}\label{le41}
		Let $m \in (m_*, m^*)$ and \eqref{gamma}. Suppose that non-negative initial data $\rho_{0\varepsilon}$ satisfies the condition \eqref{in}.
		If the initial free energy satisfies
		\begin{equation}\label{ba}
			\|\rho_{0\varepsilon}\|_{L^1(\mathbb{R}^d)}^\frac{\tau_1 m}{\tau_2} \mathcal{F_{\varepsilon}}(\rho_{0\varepsilon}) < \|V\|_{L^1(\mathbb{R}^d)}^{\frac{\tau_1m}{\tau_2}}\mathcal{F}(V)
		\end{equation}
		and the initial data satisfies
		\begin{equation}\label{bb}
			\|\rho_0\|_{L^m(\mathbb{R}^d)}<s^*,
		\end{equation}
		where $s^*$ is given in \eqref{ss}, then
		\begin{equation}\label{41}
			\rho_\varepsilon\in L^\infty\left(0,\infty;L^{p_0}(\mathbb{R}^d)\right)
		\end{equation}
		holds, where
		\begin{equation}\label{p}
			p_0=\max\{2m+2,\gamma+2, \frac{2md}{d-m(d-\gamma-2)}-1+\tau, d-\frac{dm}{2}-1+\tau\}, \quad 0<\tau\ll 1.
		\end{equation}
	\end{lemma}
	\begin{proof}
		Fristly, we provide the proof of the uniform boundedness estimate in $L^m$ for the regularized solution $\rho_{\varepsilon}$.

		According to \eqref{r2} and the monotonic decrease of the free energy, there holds
		$$
		\begin{aligned}
			\|\rho_{0\varepsilon}\|_{L^1}^\frac{\tau_1 m}{\tau_2} \mathcal{F_{\varepsilon}}(\rho_{0\varepsilon})
			\geq & \|\rho_\varepsilon\|_{L^1}^\frac{\tau_1 m}{\tau_2}\left(\int_{\mathbb{R}^d}\frac{1}{m-1}\rho_{\varepsilon}^mdx-\frac{\bar{c}_d}{2}C_{\text{HLS}}\|\rho_\varepsilon\|_{L^1}^{2(1-\theta)}\|\rho_{\varepsilon}\|_{L^m}^{2\theta}\right)\\
			=&f_2\left(\|\rho_{\varepsilon}\|_{L^1}^{\tau_1}\|\rho_{\varepsilon}\|_{L^m}^{\tau_2}\right),
		\end{aligned}
		$$
		where  $$
		f_2(s)=\frac{1}{m-1}s^{\frac{m}{\tau_2}}-\frac{1}{2}C_{\text{HLS}}\bar{c}_ds^{\frac{2\theta}{\tau_2}},
		$$
		and $\tau_1, \tau_2$ are given in \eqref{tau}.
		For all $\rho\in L^1(\mathbb{R}^d) \cap L^{m}(\mathbb{R}^d)$, it follows from \eqref{b2} that
		$$
		f_2\left(\|V\|_{L^1(\mathbb{R}^d)}^{\tau_1}\|V\|_{L^m(\mathbb{R}^d)}^{\tau_2}\right)
		\geq f_2\left(\|\rho\|_{L^1(\mathbb{R}^d)}^{\tau_1} \|\rho\|_{L^m(\mathbb{R}^d)}^{\tau_2}\right).
		$$
		It is obvious that $f_2(x)$ reaches its maximum value at
		\begin{equation}\label{b3}
			\tilde{s}_2=\left(\frac{m}{\theta C_{\text{HLS}}\bar{c}_d(m-1)}\right)^{\frac{\tau_2}{2\theta-m}}=\left(\frac{2d}{\gamma C_{\text{HLS}}\bar{c}_d}\right)^{\frac{\tau_1}{2(1-\theta)}}=\|V\|_{L^1(\mathbb{R}^d)}^{\tau_1}\|V\|_{L^m(\mathbb{R}^d)}^{\tau_2}.
		\end{equation}
		Thus, combining the results above with \eqref{ba}, we have
		$$
		f_2(\tilde{s}_2)> f_2\left(\|\rho_{\varepsilon}\|_{L^1}^{\tau_1}\|\rho_{\varepsilon}\|_{L^m}^{\tau_2}\right).
		$$
		Since $f_2(x)$ is a monotone increasing function over $0<x<\tilde{s}_2$ and \eqref{bb} holds, we can take inverse function to obtain the absolute constant $\mu_3<1$ such that
		\begin{equation}\label{b9}
			\|\rho_{\varepsilon}\|_{L^1}^{\tau_1}\|\rho_{\varepsilon}\|_{L^m}^{\tau_2}
			<\mu_3\tilde{s}_2, \quad \text{for all } t\geq 0.
		\end{equation}

		Next, we prove the uniform $L^{p_0}$ bound for the regularized  solution $\rho_\varepsilon$.
		
		Taking $p_0\rho_{\varepsilon}^{p_0-1}$ with \eqref{p} as a test function in \eqref{r1} and integrating the result equation in $\mathbb{R}^d$, we have
		\begin{equation}\label{b5}
			\begin{aligned}
				&\frac{d}{dt}\int_{\mathbb{R}^d}\rho_\varepsilon^{p_0}dx+\frac{4pm(p_0-1)}{(m+p_0-1)^2}\int_{\mathbb{R}^d}|\nabla \rho_\varepsilon^{\frac{m+p_0-1}{2}}|^2dx\\
				\leq &(p_0-1)\int_{\mathbb{R}^d}\rho_{\varepsilon}^{p_0}(-\Delta c_\varepsilon)dx
				+\varepsilon p_0\int_{\mathbb{R}^d}\rho_\varepsilon^{p_0-1}(-\Delta c_\varepsilon)dx:=K_1+K_2.
			\end{aligned}
		\end{equation}
		According to \eqref{r28}, it follows from \eqref{hls} that
		\begin{equation}\label{b4}
			K_1\leq (p_0-1)d\gamma \bar{c}_dC_{\text{HLS}}\|\rho_\varepsilon^{p_0}\|_{L^{r_3}}\|\rho_\varepsilon\|_{L^{s_3}},
		\end{equation}
		where $\frac{1}{r_3}+\frac{1}{s_3}+\frac{\gamma+2}{d}=2$.
		Take $r_3=\frac{p_0+1}{p_0}$. Then
		we have from \eqref{p} that $m<s_3=\frac{d(p_0+1)}{(d-\gamma-2)(p_0+1)+d}<p_0+1$. Using the interpolation inequality, we deduce that
		$$
		K_1\leq (p_0-1)d\gamma\bar{c}_dC_{\text{HLS}}\|\rho_\varepsilon\|_{L^{p_0+1}}^{p_0+\lambda_1}\|\rho_\varepsilon\|_{L^m}^{1-\lambda_1}, \quad \lambda_1=\frac{(s_3-m)(p_0+1)}{s_3(p_0+1-m)}.
		$$
		Similarly, we can obtain from $r_4=\frac{p_0+1}{p_0-1}$ and $s_4=\frac{d(p_0+1)}{(d-\gamma-2)(p_0+1)+2d}$ that
		\begin{equation}\label{b6}
			\begin{aligned}
				K_2\leq &\varepsilon p_0d\gamma\bar{c}_d C_{\text{HLS}}\|\rho_\varepsilon^{p_0-1}\|_{L^{r_4}}\|\rho_\varepsilon\|_{L^{s_4}}\\
				\leq & \varepsilon p_0d\gamma\bar{c}_d C_{\text{HLS}}\|\rho_\varepsilon\|_{L^{p_0+1}}^{p_0-1+\lambda_2}\|\rho_\varepsilon\|_{L^m}^{1-\lambda_2}, \quad \lambda_2=\frac{(s_4-m)(p_0+1)}{s_4(p_0+1-m)}.
			\end{aligned}
		\end{equation}
		Since \eqref{p} leads to $m<p_0+1<\frac{d(m+p_0-1)}{d-2}$, we have
		$$
		\|\rho_\varepsilon\|_{L^{p_0+1}}\leq \|\rho_\varepsilon\|_{L^m}^{\lambda_3}\left(S_d^{-\frac{1}{2}}\|\nabla \rho_\varepsilon^{\frac{m+p_0-1}{2}}\|_{L^2}\right)^{\frac{2(1-\lambda_3)}{m+p_0-1}}, \quad \lambda_3=\frac{m[2(p_0+1)-d(2-m)]}{(p_0+1)[d(p_0-1)+2m]}.
		$$
		Clearly, $p_0-1+\lambda_2<p_0+\lambda_1$ and $m>\frac{2d}{2d-\gamma}$ yield
		$$
		\frac{2(p_0-1+\lambda_2)(1-\lambda_3)}{m+p_0-1}
		<\frac{2(p_0+\lambda_1)(1-\lambda_3)}{m+p_0-1}<2,
		$$
		which implies from Young's inequality that
		\begin{equation}\label{b7}
			\begin{aligned}
				&(p_0-1)d\gamma\bar{c}_dC_{\text{HLS}}\|\rho_\varepsilon\|_{L^{p_0+1}}^{p_0+\lambda_1}\|\rho_\varepsilon\|_{L^m}^{1-\lambda_1}\\
				\leq & \frac{p_0m(p_0-1)}{(m+p_0-1)^2}\|\nabla \rho_\varepsilon^{\frac{m+p_0-1}{2}}\|_{L^2}^{2}
				+C\left(\gamma, d, m, p_0, \mu_3, \tilde{s}_2, \|\rho_0\|_{L^1}\right)
			\end{aligned}
		\end{equation}
		and
		\begin{equation}\label{b8}
			\begin{aligned}
				&\varepsilon p_0\bar{c}_dd\gamma C_{\text{HLS}}\|\rho_\varepsilon\|_{L^{p_0+1}}^{p_0-1+\lambda_2}\|\rho_\varepsilon\|_{L^m}^{1-\lambda_2}\\
				\leq &\frac{p_0m(p_0-1)}{(m+p_0-1)^2}\|\nabla \rho_\varepsilon^{\frac{m+p_0-1}{2}}\|_{L^2}^{2}
				+C\left(\gamma, d, m, p_0, \mu_3, \tilde{s}_2, \|\rho_0\|_{L^1}\right).
			\end{aligned}
		\end{equation}
		Combining \eqref{b4}-\eqref{b8} and \eqref{b9} with \eqref{b5}, we get
		\begin{equation}\label{b10}
			\frac{d}{dt}\int_{\mathbb{R}^d}\rho_\varepsilon^{p_0}dx
			\leq -\frac{2p_0m(p_0-1)}{(m+p_0-1)^2}\int_{\mathbb{R}^d}\big|\nabla \rho_\varepsilon^{\frac{m+p_0-1}{2}}\big|^2dx+C\left(\gamma, d, m, p_0, \mu_3, \tilde{s}_2, \|\rho_0\|_{L^1}\right).
		\end{equation}
		Applying interpolation inequality with $1<p_0<\frac{(m+p_0-1)d}{d-2}$, \eqref{So} and Young's inequality with $\frac{2p_0(1-\lambda_4)}{m+p_0-1}<2$, we obtain
		$$
		\begin{aligned}
			\|\rho_\varepsilon\|_{L^{p_0}}^{p_0}
			\leq & \|\rho_\varepsilon\|_{L^1}^{p_0\lambda_4}\left(S_d^{-\frac{1}{2}}\|\nabla \rho_\varepsilon^{\frac{m+p_0-1}{2}}\|_{L^2}\right)^{\frac{2p_0(1-\lambda_4)}{m+p_0-1}}\\
			\leq &\frac{2p_0m(p_0-1)}{(m+p_0-1)^2}\|\nabla \rho_\varepsilon^{\frac{m+p_0-1}{2}}\|_{L^2}^{2}
			+C\left(d,m,p_0, \|\rho_0\|_{L^1}\right).
		\end{aligned}
		$$
		Here $\lambda_4=\frac{d(m-1)+2_0}{p_0[d(m+p_0-2)+2]}$.
		Inserting the result above into \eqref{b10}, we obtain
		$$
		\frac{d}{dt}\|\rho_\varepsilon\|_{L^{p_0}}^{p_0}
		+\|\rho_\varepsilon\|_{L^{p_0}}^{p_0}
		\leq  C\left(\gamma, d, m, p_0, \mu_3, \tilde{s}_2, \|\rho_0\|_{L^1}\right).
		$$
		Thus, taking integral with respect to that time variable, and using \eqref{b9}, we derive that
		$$
		\|\rho_\varepsilon\|_{L^{p_0}}^{p_0}\leq \|\rho_{0\varepsilon}\|_{L^{p_0}}^{p_0}+C\left(\gamma, d, m, p_0, \mu_3, \tilde{s}_2, \|\rho_0\|_{L^1}\right), \quad \mbox{for all } t\geq 0.
		$$
		Therefore,  we can easily obtain \eqref{41}.
	\end{proof}		
	
	We now provide additional moment estimates for $\rho_\varepsilon^m$, which are crucial for compactness arguments.
	\begin{lemma}\label{le42}
		Under the assumptions of \eqref{in}, we have for any $t \in [0,T]$ with $T>0$,
		\begin{equation}\label{b41}
			\int_{\mathbb{R}^d}(1+|x|)\rho_{\varepsilon}^{m}(x,t)dx \leq C,
		\end{equation}
		where $C$ depends on $d,\gamma,m, \mu_3, \tilde{s}_2, \|\rho_0\|_{L^1(\mathbb{R}^d)}, f_2(\tilde{s}_2),  \|\rho_{\varepsilon}\|_{L^{p_0}(\mathbb{R}^d)}, T$ and $\int_{\mathbb{R}^d}|x|^2\rho_0(x)dx$.
	\end{lemma}
	
	\begin{proof}
		Since $\int_{\mathbb{R}^d}\rho_{\varepsilon}^{m}(x,t)dx\leq C\left(\mu_3, x_0, \|\rho_0\|_{L^1}\right)$ from \eqref{b9}, we focus on proving $$
\int_{\mathbb{R}^d}|x|\rho_{\varepsilon}^{m}(x,t)dx \leq C, \quad {\rm for \ any} \ \ t\in [0,T].
$$
		
		Applying the H\"{o}lder inequality, we have
		\begin{equation}\label{b79}
			\int_{\mathbb{R}^{d}} |x|\rho_{\varepsilon}^{m} dx
			\leq \|\rho_\varepsilon\|_{L^{2m-1}}^{m-\frac{1}{2}}\left(\int_{\mathbb{R}^{d}} |x|^{2} \rho_{\varepsilon}dx\right)^{\frac{1}{2}}.
		\end{equation}
		From \eqref{r1} and using integration by parts, we obtain
		\begin{equation}\label{b77}
			\begin{aligned}
				\int_{\mathbb{R}^{d}} |x|^{2} \partial_{t} \rho_{\varepsilon} dx
				&= \int_{\mathbb{R}^{d}} \Delta |x|^{2} \left[ (\rho_{\varepsilon} + \varepsilon)^{m} - \rho_{\varepsilon}^{m} \right] dx + \int_{\mathbb{R}^{d}} \nabla |x|^{2} \cdot (\rho_{\varepsilon} + \varepsilon) \nabla c_{\varepsilon} dx \\
				&=2d \int_{\mathbb{R}^{d}} \left[ (\rho_{\varepsilon} + \varepsilon)^{m} - \rho_{\varepsilon}^{m} \right] dx - 2d \varepsilon \int_{\mathbb{R}^{d}} c_{\varepsilon}  dx + 2 \int_{\mathbb{R}^{d}} x \cdot \nabla c_{\varepsilon} \rho_{\varepsilon}  dx.
			\end{aligned}
		\end{equation}
		The fact $(x-y)\cdot x = \frac{|x|^{2} - |y|^{2} + |x - y|^{2}}{2}$, together with \eqref{hls2}, yields
		\begin{equation}\label{t}
			\begin{aligned}
				&2 \int_{\mathbb{R}^{d}} x \cdot \nabla c_{\varepsilon} \rho_{\varepsilon} \, dx \\
				&=-\bar{c}_{d}\gamma\int_{\mathbb{R}^{d}} \int_{\mathbb{R}^{d}} \frac{(x-y)\cdot x \rho_{\varepsilon}(x) \rho_{\varepsilon}(y)}{(|x - y|^{2} + \varepsilon^{2})^{\frac{\gamma + 2}{2}}} dx dy \\
				&= -\bar{c}_{d} \gamma \int_{\mathbb{R}^{d}} \int_{\mathbb{R}^{d}} \frac{(|x - y|^{2} + \varepsilon^{2}) \rho_{\varepsilon}(x) \rho_{\varepsilon}(y)}{(|x - y|^{2} + \varepsilon^{2})^{\frac{\gamma + 2}{2}}} dx dy + \bar{c}_{d} \gamma \varepsilon^{2} \int_{\mathbb{R}^{d}} \int_{\mathbb{R}^{d}} \frac{\rho_{\varepsilon}(x) \rho_{\varepsilon}(y)}{(|x-y|^{2} + \varepsilon^{2})^{\frac{\gamma + 2}{2}}} dx dy \\
				&\leq-\gamma\int_{\mathbb{R}^{d}}\rho_{\varepsilon} c_{\varepsilon} dx +\bar{c}_{d}\gamma \varepsilon^{2} \|\rho_{\varepsilon}\|_{L^{\frac{2d}{2d-\gamma-2}}}^{2}.
			\end{aligned}
		\end{equation}
		By the interpolation inequality with $1<\frac{2d}{2d-\gamma-2}<p_0$, we have
		$$
		\|\rho_\varepsilon\|_{L^{\frac{2d}{2d-\gamma-2}}}\leq C\left(d, \gamma, m, \|\rho_0\|_{L^1},\|\rho_\varepsilon\|_{L^{p_0}}\right).
		$$
		Inserting this result, \eqref{r2} and \eqref{t} into \eqref{b77}, and using \eqref{41}, we obtain
		$$\begin{aligned}
		\int_{\mathbb{R}^{d}} |x|^{2} \partial_{t} \rho_{\varepsilon}  dx 
&\leq \left(2d - \frac{2\gamma}{m - 1}\right) \int_{\mathbb{R}^{d}} 
\left[ (\rho_{\varepsilon} + \varepsilon)^{m} - \rho_{\varepsilon}^{m} \right] dx \\
&+ 2\gamma \mathcal{F}_{\varepsilon}(\rho_{\varepsilon})
		+ C\left(d,\gamma, m,  \|\rho_0\|_{L^1}, \|\rho_\varepsilon\|_{L^{p_0}}\right).
\end{aligned}		$$
		Clearly, the fact $m < \frac{d + \gamma}{d}$ leads to $2d - \frac{2\gamma}{m - 1} < 0$. From the monotonicity of free energy and \eqref{ba}, there holds
		$$
		\int_{\mathbb{R}^{d}} |x|^{2} \partial_{t} \rho_{\varepsilon} dx
		\leq 2\gamma \mathcal{F}_\varepsilon(\rho_{0\varepsilon})+C\left(d,\gamma,\|\rho_0\|_{L^1(\mathbb{R}^d)}, \|\rho_\varepsilon\|_{L^{p_0}}\right)
		\leq C\left(d, \gamma, m, f_2(\tilde{s}_2), \|\rho_0\|_{L^1},  \|\rho_{\varepsilon}\|_{L^{p_0}}\right).
		$$
		which implies
		\begin{equation}\label{r58}
			\int_{\mathbb{R}^d}|x|^{2} \rho_{\varepsilon}dx \leq \int_{\mathbb{R}^d}|x|^2\rho_{0\varepsilon}dx+C\left(d,\gamma,m,   f_2(\tilde{s}_2), \|\rho_0\|_{L^1},  \|\rho_{\varepsilon}\|_{L^{p_0}}, T\right).
		\end{equation}
		Applying the interpolation inequality to $\|\rho_\varepsilon\|_{L^{2m-1}}$ with $1<2m-1<p_0$, we have
		$$
		\|\rho_\varepsilon\|_{L^{2m-1}}\leq \|\rho_\varepsilon\|_{L^1}^{\theta_5}\|\rho_\varepsilon\|_{L^{p_0}}^{1-\theta_5}, \quad \theta_5=\frac{p_0-(2m-1)}{(2m-1)(p_0-1)}.
		$$
		Combining this result, \eqref{41} and \eqref{r58} with \eqref{b79}, we obtain that $\int_{\mathbb{R}^d}|x|\rho_\varepsilon^m(x,t)dx\leq C$ for any $t\in [0,T]$. Thus, we complete the proof of Lemma \ref{le42}.
	\end{proof}

	With the key lemmas above, we now proceed to the proof of the global existence part of Theorem \ref{th2}.

	\emph{Proof of Theorem \ref{th2} (1).}
	
	According to Lemma \ref{le41}, there exists a subsequence $\rho_{\varepsilon}$, not relabeled, such that
	$$
	\rho_\varepsilon \rightharpoonup \rho \quad \text{weak star in } L^\infty\left(0,T; L^p(\mathbb{R}^d)\right).
	$$
	We claim that there exists a subsequences of $\rho_\varepsilon$ denoted again by $\rho_\varepsilon$ such that for $\frac{2}{m}<\bar{p}''<2$
	$$
	\rho_\varepsilon \to \rho \quad \text{strongly in } C\left(0,T; L^{m\bar{p}''}(\mathbb{R}^d)\right),
	$$
	$$
	\nabla \rho_\varepsilon^m \rightharpoonup  \nabla \rho^m  \quad \text{weak star in } L^\infty\left(0,T; L^2(\mathbb{R}^d)\right).
	$$

	In fact, multiplying $\partial_t(\rho_\varepsilon+\varepsilon)^m$ on both sides of \eqref{b1}, and integrating with respect to the space variable over $\mathbb{R}^d$, we obtain
	\begin{equation}\label{i}
		\begin{aligned}
			&\frac{4m}{(m+1)^2}\int_{\mathbb{R}^d}\left|\partial_t(\rho_\varepsilon+\varepsilon)^{\frac{m+1}{2}}\right|^2dx
			+\frac{1}{2}\frac{d}{dt}\int_{\mathbb{R}^d}\left|\nabla (\rho_\varepsilon+\varepsilon)^m\right|^2dx\\
			= & -\int_{\mathbb{R}^d}\nabla (\rho_\varepsilon+\varepsilon)\cdot \nabla c_\varepsilon \partial_t(\rho_\varepsilon+\varepsilon)^mdx
			-\int_{\mathbb{R}^d}(\rho_\varepsilon+\varepsilon)\Delta c_\varepsilon\partial_t(\rho_\varepsilon+\varepsilon)^mdx:=I_1+I_2.
		\end{aligned}
	\end{equation}
	According to \eqref{gamma} and \eqref{41}, we deduce from Young's inequality that
	$$
	\begin{aligned}
		\|\nabla c_\varepsilon\|_{L^{\infty}}
		\leq & \bar{c}_d\gamma
		\left \| \int_{|x-y|\leq 1} \frac{\rho_\varepsilon(x)}{|x-y|^{\gamma+1}}dx
		+\int_{|x-y|> 1} \frac{\rho_\varepsilon(x)}{|x-y|^{\gamma+1}}dx \right \|_{L^\infty}\\
		\leq &  \bar{c}_d\gamma
		\left \|\frac{\gamma+1}{\gamma+2}\int_{|x-y|\leq 1}\frac{1}{|x-y|^{\gamma+2}}dx
		+\frac{1}{\gamma+2}\int_{\mathbb{R}^d}\rho_\varepsilon^{\gamma+2}dx+\|\rho_\varepsilon\|_{L^1}\right \|_{L^\infty}.
	\end{aligned}
	$$
	If $\gamma+2=p_0$, we have $\|\nabla c_\varepsilon\|_{L^{\infty}}\leq C\left(d,\gamma, \|\rho_\varepsilon\|_{L^1}, \|\rho_\varepsilon\|_{L^{p_0}} \right)$; otherwise, by the interpolation inequality with $1<\gamma+2<p_0$, we also have $\|\nabla c_\varepsilon\|_{L^{\infty}}\leq C\left(d,\gamma, \|\rho_\varepsilon\|_{L^1}, \|\rho_\varepsilon\|_{L^{p_0}} \right)$.
	Thus, using Young's inequality again, we obtain that
	\begin{equation}\label{i1}
		\begin{aligned}
			|I_1|\leq & \frac{m}{2}\|\nabla c_\varepsilon\|_{L^\infty}^2\int_{\mathbb{R}^d} (\rho_\varepsilon+\varepsilon)^{m-1}|\nabla(\rho_\varepsilon+\varepsilon)|^2dx+\frac{m}{2}\int_{\mathbb{R}^d}(\rho_\varepsilon+\varepsilon)^{m-1}|\partial_t(\rho_\varepsilon+\varepsilon)|^2dx\\
			\leq & \frac{2m}{(m+1)^2}\|\nabla c_\varepsilon\|_{L^\infty}^2 \int_{\mathbb{R}^d}|\nabla(\rho_\varepsilon+\varepsilon)^{\frac{m+1}{2}}|^2dx+\frac{2m}{(m+1)^2}\int_{\mathbb{R}^d}\left|\partial_t(\rho_\varepsilon+\varepsilon)^{\frac{m+1}{2}}\right|^2dx.
		\end{aligned}
	\end{equation}
	According to \eqref{gamma} and \eqref{41}, we deduce from the weak Young inequality \cite[Page 107]{LL} and the interpolation inequality with $1<p',p''<p_0$ that
	$$
	\|\Delta c_\varepsilon\|_{L^2}\leq \bar{c}_dd\gamma\left\| \frac{1}{|x|^{\gamma+2}}\ast \rho_\varepsilon\right\|_{L^2}
	\leq \bar{c}_dd\gamma \left\|\frac{1}{|x|^{\gamma+2}}\right\|_{L_w^{\frac{d}{\gamma+2}}}
	\|\rho_\varepsilon\|_{L^{p'}}
	\leq C\left(d, \gamma, \|\rho_0\|_{L^1}, \|\rho_\varepsilon\|_{L^{p_0}}\right)
	$$
	and
	$$
	\|\Delta c_\varepsilon\|_{L^4}\leq \bar{c}_dd\gamma\left\| \frac{1}{|x|^{\gamma+2}}\ast \rho_\varepsilon\right\|_{L^4}
	\leq \bar{c}_dd\gamma \left\|\frac{1}{|x|^{\gamma+2}}\right\|_{L_w^{\frac{d}{\gamma+2}}}
	\|\rho_\varepsilon\|_{L^{p''}}
	\leq C\left(d, \gamma, \|\rho_0\|_{L^1}, \|\rho_\varepsilon\|_{L^{p_0}}\right).
	$$
	Here $p'=\frac{2d}{3d-2\gamma-4}>1$ and $p''=\frac{4d}{5d-4(\gamma+2)}>1$. Thus, we derive from Young's inequality with $0<\delta<1$ that
	\begin{equation}\label{i2}
		|I_2|\leq  \frac{2m\delta}{(m+1)^2}\int_{\mathbb{R}^d}\left|\partial_t(\rho_\varepsilon+\varepsilon)^{\frac{m+1}{2}}\right|^2dx
		+\frac{m}{2\delta}\int_{\mathbb{R}^d}(\rho_\varepsilon+\varepsilon)^{m+1}|\Delta c_\varepsilon|^2dx,
	\end{equation}
	where from \eqref{41}, we have
	\begin{equation}\label{i22}
		\begin{aligned}
			&\frac{m}{2\delta}\int_{\mathbb{R}^d}(\rho_\varepsilon+\varepsilon)^{m+1}|\Delta c_\varepsilon|^2dx\\
			\leq & \frac{m}{2\delta}2^{m-1}\int_{\mathbb{R}^d}\rho_\varepsilon^{m+1}|\Delta c_\varepsilon|^2dx
			+\frac{m}{2\delta}2^{m-1}\varepsilon^{m+1}\int_{\mathbb{R}^d}|\Delta c_\varepsilon|^2dx\\
			\leq & \frac{m}{2^{3-m}\delta}\int_{\mathbb{R}^d}\rho_\varepsilon^{2m+2}dx
			+\frac{m}{2^{3-m}\delta}\int_{\mathbb{R}^d}|\Delta c_\varepsilon|^4dx
			+C\left(d,\gamma, m, \|\rho_0\|_{L^1}, \|\rho_\varepsilon\|_{L^{p_0}}\right).
		\end{aligned}
	\end{equation}
	If $2m+2=p_0$, we have $\frac{m}{2\delta}\int_{\mathbb{R}^d}(\rho_\varepsilon+\varepsilon)^{m+1}|\Delta c_\varepsilon|^2dx\leq C\left(d,\gamma, m, \|\rho_0\|_{L^1}, \|\rho_\varepsilon\|_{L^{p_0}} \right)$; otherwise, by the interpolation inequality with $1<2m+2<p_0$, we can also obtain this result.
	Combining this estimate and \eqref{i1}-\eqref{i22} with \eqref{i}, we obtain
	$$
	\begin{aligned}
		&\frac{2m(1-\delta)}{(m+1)^2}\int_{\mathbb{R}^d}\left|\partial_t(\rho_\varepsilon+\varepsilon)^{\frac{m+1}{2}}\right|^2dx
		+\frac{1}{2}\frac{d}{dt}\int_{\mathbb{R}^d}\left|\nabla (\rho_\varepsilon+\varepsilon)^m\right|^2dx\\
		\leq & \frac{2m}{(m+1)^2}\|\nabla c_\varepsilon\|_{L^{\infty}}^2\int_{\mathbb{R}^d}\left|\nabla(\rho_\varepsilon+\varepsilon)^{\frac{m+1}{2}}\right|^2dx
		+C\left(d,\gamma, m, \|\rho_0\|_{L^1}, \|\rho_\varepsilon\|_{L^{p_0}} \right).
	\end{aligned}
	$$
	Integrating the inequality above with respect to time variable, and noticing \eqref{41}, we get
	\begin{equation}\label{b53}
		\begin{aligned}
			&\frac{2m(1-\delta)}{(m+1)^2}\int_0^T\int_{\mathbb{R}^d}\left|\partial_t(\rho_\varepsilon+\varepsilon)^{\frac{m+1}{2}}\right|^2dx
			+\frac{1}{2}\sup_{0<t<T}\int_{\mathbb{R}^d}\left|\nabla (\rho_\varepsilon+\varepsilon)^m\right|^2dx\\
			\leq &
			\frac{1}{2}\int_{\mathbb{R}^d}\left|\nabla (\rho_{0\varepsilon}+\varepsilon)^m\right|^2dx
			+\frac{2m}{(m+1)^2}\|\nabla c_\varepsilon\|_{L^{\infty}}^2\int_0^T\int_{\mathbb{R}^d}\left|\nabla(\rho_\varepsilon+\varepsilon)^{\frac{m+1}{2}}\right|^2dx\\
			&+C\left(d,\gamma, m, \|\rho_0\|_{L^1}, \|\rho_\varepsilon\|_{L^{p_0}} \right).
		\end{aligned}
	\end{equation}

	Next, multiplying $\rho_\varepsilon$ on both sides of \eqref{b1} simultaneously, and integrating with respect to $x$ and $t$, we have
	\begin{equation}\label{b54}
		\begin{aligned}
			&\frac{4m}{(m+1)^2}\int_0^T\int_{\mathbb{R}^d}\left|\nabla(\rho_\varepsilon+\varepsilon)^{\frac{m+1}{2}}\right|^2dxdt\\
			\leq & \frac{1}{2}\|\rho_{0\varepsilon}\|_{L^2}^2
			+\frac{1}{4}\int_0^T\int_{\mathbb{R}^d}\rho_\varepsilon^4dxdt
			+\frac{1}{4}\int_0^T\int_{\mathbb{R}^d}|\Delta c_\varepsilon|^2dxdt
			+\frac{\varepsilon}{2}\int_0^T\int_{\mathbb{R}^d}\left(\rho_\varepsilon^2+|\Delta c_\varepsilon|^2\right)dxdt.
		\end{aligned}
	\end{equation}
	By the interpolation inequality with $1<4<p_0$, we derive that
	$$
	\frac{4m}{(m+1)^2}\int_0^T\int_{\mathbb{R}^d}\left|\nabla(\rho_\varepsilon+\varepsilon)^{\frac{m+1}{2}}\right|^2dxdt
	\leq C\left(d,\gamma, m, \|\rho_0\|_{L^1}, \|\rho_\varepsilon\|_{L^{p_0}}, \|\rho_{0 \varepsilon}\|_{L^2} \right).
	$$
	Inserting \eqref{b54} and \eqref{na} into \eqref{b53}, and using \eqref{41}, we have
	\begin{equation}\label{k1}
		\rho_\varepsilon^m\in L^\infty{\left(0,\infty;H^1(\mathbb{R}^d)\right)},
	\end{equation}
	and
	\begin{equation}\label{k2}
		\partial_t(\rho_\varepsilon+\varepsilon)^{\frac{m+1}{2}}\in L^2\left(0,\infty;L^2(\mathbb{R}^d)\right).
	\end{equation}
	The fact
	$\partial_t\rho_\varepsilon^m
	\leq \frac{2m}{m+1}\rho_\varepsilon ^{\frac{m-1}{2}} \partial_t(\rho_\varepsilon+\varepsilon)^{\frac{m+1}{2}}$
	yields that for $\frac{2}{m}<\bar{p}''<2$
	\begin{equation}\label{b55}
		\begin{aligned}
			&\int_0^T\int_{\mathbb{R}^d}\left|\partial_t\rho_\varepsilon^m\right|^{\bar{p}''}dxdt\\
			\leq &
			\left(\frac{2m}{m+1}\right)^{\bar{p}''}\left(\frac{\bar{p}''}{2}
			\int_0^T\int_{\mathbb{R}^d}\left| \partial_t(\rho_\varepsilon+\varepsilon)^{\frac{m+1}{2}}\right|^2dxdt
			+\frac{2-\bar{p}''}{2}\int_{\mathbb{R}^d}\rho_\varepsilon^{\frac{(m-1)\bar{p}''}{2-\bar{p}''}}dx\right).
		\end{aligned}
	\end{equation}
	By the interpolating inequality with $1<\bar{p}''<p_0$, and \eqref{k2}, we deduce that 
$$
\partial_t\rho_\varepsilon^m\in L^{\bar{p}''}\left(0,\infty;L^{\bar{p}''}(\mathbb{R}^d)\right).
$$ 
According to  \eqref{k1}, \cite[Lemma 5.7]{CWW} and \cite[Lemma 3.4]{CNP}, we have
	$$
	H^1(\mathbb{R}^d)\cap L^1\left(\mathbb{R}^d; (1+|x|)dx\right) \hookrightarrow\hookrightarrow L^{\bar{p}''}(\mathbb{R}^d)\hookrightarrow L^{\bar{p}''}(\mathbb{R}^d),
	$$
	so we satisfy the conditions of Lions-Aubin Lemma \cite{A} (see also \cite[Chapter IV, Section 4]{Lions}), which implies that there exists a subsequence of $\rho_\varepsilon$, not relabeled, such that
	$$
	\rho_\varepsilon^m \to \varpi \quad \text{strongly in } C\left(0,\infty; L^{\bar{p}''}(\mathbb{R}^d)\right).
	$$
	Taking $\varpi=\rho^m$, together with the fact $|b-a|^m\leq |b^m-a^m|$ for $0\leq a \leq b$ and $m>1$, we deduce
	$$
	\rho_\varepsilon \to \rho \quad \text{strongly in } C\left(0,\infty; L^{m\bar{p}''}(\mathbb{R}^d)\right)
	$$
	and
	$$
	\nabla \rho_\varepsilon^m \rightharpoonup  \nabla \rho^m  \quad \text{weak star in } L^\infty\left(0,\infty; L^2(\mathbb{R}^d)\right).
	$$
	Thus, the claim holds. Similar to the argument in \cite{BL,CLW}, we take limit of \eqref{r1} and obtain that for any test function $\psi\in C_0^\infty$ in Definition \ref{de}, there exists a global weak solution $\rho$ of \eqref{a1}. Thus, we complete the proof of Theorem \ref{th2}.
	\qed
	
	\section{The equivalence of sharp criteria.}
In this section, we discuss the equivalence of the sharp criteria presented in Theorem \ref{th1} and Theorem \ref{th2} regarding global existence and finite-time blow-up of solutions.

	\begin{proof}
		Firstly, we prove the equivalence of the conditions \eqref{c1} and \eqref{c2}.
		
		It is enough to show that 	
		\begin{equation}\label{conditions}
			\frac{\gamma-d(m-1)}{\gamma(m-1)} \left(\frac{2d^2\alpha(d)}{ C_{\text{HLS}}} \right)^{\frac{d(m-1)}{\gamma-d(m-1)}} \|\rho_0\|_{L^1(\mathbb{R}^d)}^{\frac{2d-m(2d-\gamma)}{\gamma-d(m-1)}}=	\|\rho_0\|_{L^1(\mathbb{R}^d)}^{-\frac{\tau_1 m}{\tau_2}} \|V\|_{L^1(\mathbb{R}^d)}^{\frac{\tau_1m}{\tau_2}}\mathcal{F}(V).
		\end{equation}
		According to \eqref{F} and \eqref{b2}, we have
		$$
		\|V\|_{L^1(\mathbb{R}^d)}^{\frac{\tau_1m}{\tau_2}}\mathcal{F}(V)
		 = \frac{1}{m-1}\left(\|V\|_{L^1(\mathbb{R}^d)}^{\tau_1}\|V\|_{L^m(\mathbb{R}^d)}^{\tau_2}\right)^{\frac{m}{\tau_2}}-\frac{1}{2}C_{\text{HLS}}\bar{c}_d\left(\|V\|_{L^1(\mathbb{R}^d)}^{\tau_1}\|V\|_{L^m(\mathbb{R}^d)}^{\tau_2}\right)^{\frac{2\theta}{\tau_2}}.
		$$
		Thus, we derive from the definition of $\bar{c}_d$, \eqref{tau} and \eqref{b3} that
		$$
		\begin{aligned}
			&\|\rho_0\|_{L^1(\mathbb{R}^d)}^{-\frac{\tau_1 m}{\tau_2}} \|V\|_{L^1(\mathbb{R}^d)}^{\frac{\tau_1m}{\tau_2}}\mathcal{F}(V)\\
			=& \|\rho_0\|_{L^1(\mathbb{R}^d)}^{\frac{2d-m(2d-\gamma)}{\gamma-d(m-1)}}
			\left[\frac{1}{m-1}\left(\frac{2d}{\gamma C_{\text{HLS}}\bar{c}_d}\right)^{\frac{\tau_1m}{2\tau_2(1-\theta)}}
			-\frac{1}{2}C_{\text{HLS}}\bar{c}_d\left(\frac{2d}{\gamma C_{\text{HLS}}\bar{c}_d}\right)^{\frac{2\theta}{2\theta-m}}\right]\\
			=& \|\rho_0\|_{L^1(\mathbb{R}^d)}^{\frac{2d-m(2d-\gamma)}{\gamma-d(m-1)}}
			\left[\frac{1}{m-1}\left(\frac{2d}{\gamma C_{\text{HLS}}\bar{c}_d}\right)^{\frac{d(m-1)}{\gamma-d(m-1)}}-
			\frac{d}{\gamma}\left(\frac{2d}{\gamma C_{\text{HLS}}\bar{c}_d}\right)^{\frac{\gamma}{\gamma-d(m-1)}-1}\right],
		\end{aligned}
		$$
		which implies that \eqref{conditions} holds.
		
		Next, we show that both criteria for global existence and blow-up of solutions are equivalent.
		
		If the condition for global existence of solutions in Theorem \ref{th2} holds, i.e.,
		\begin{equation}\label{2}
			\|\rho_0\|_{L^m(\mathbb{R}^d)} < s^*
			= \left(\frac{2d^2\alpha(d)}{ C_{\text{HLS}}}\right)^{\frac{d(m-1)}{\gamma m-md(m-1)}}\|\rho_0\|_{L^1(\mathbb{R}^d)}^{\frac{\gamma m-2d(m-1)}{\gamma m-md(m-1)}},
		\end{equation}
		we obtain from \eqref{int} and \eqref{2} that
		\begin{equation}\label{1}
			\|\rho_0\|_{L^{m_*}(\mathbb{R}^d)}
			\leq \|\rho_0\|_{L^1}^{1-\theta}\|\rho_0\|_{L^m}^{\theta}
			< \left( \frac{2d^2 \alpha(d)}{C_{\text{HLS}}}\right)^{\frac{\gamma}{2\gamma-2d(m-1)}}
			\|\rho_0\|_{L^1(\mathbb{R}^d)}^{\frac{\gamma m-2d(m-1)}{2\gamma-2d(m-1)}}=s_*.
		\end{equation}
		This means that global existence of solutions in Theorem \ref{th1} can be obtained.
		
		On the other hand, if the condition for finite-time blow-up of solutions in Theorem \ref{th1} holds, i.e.,
		\begin{equation}\label{3}
			\|\rho_0\|_{L^{m_*}(\mathbb{R}^d)} > s_*
			=\left( \frac{2d^2 \alpha(d)}{C_{\text{HLS}}}\right)^{\frac{\gamma}{2\gamma-2d(m-1)}}
			\|\rho_0\|_{L^1(\mathbb{R}^d)}^{\frac{\gamma m-2d(m-1)}{2\gamma-2d(m-1)}},
		\end{equation}
		we obtain from \eqref{int} and \eqref{3} that
		\begin{equation}\label{4}
			\|\rho_0\|_{L^m(\mathbb{R}^d)}
			\geq \|\rho_0\|_{L^1}^{\frac{\theta-1}{\theta}}
			\|\rho_0\|_{L^{m_*}}^{\frac{1}{\theta}}
			>
			\left(\frac{2d^2\alpha(d)}{ C_{\text{HLS}}}\right)^{\frac{d(m-1)}{\gamma m-md(m-1)}}\|\rho_0\|_{L^1(\mathbb{R}^d)}^{\frac{\gamma m-2d(m-1)}{\gamma m-md(m-1)}}=s^*,
		\end{equation}
    	This means that the finite-time blow-up of solutions in Theorem \ref{th2} can be obtained.
    	Furthermore in both criteria ``$=$'' case do not happen under condition \eqref{c1}(or \eqref{c2}). In other words, the above criteria are both sharp. Based on \eqref{int}, any initial data fulfilling criterion \eqref{2} necessarily comply with \eqref{1}, thereby ensuring the global existence of solutions. Conversely, from the blow-up perspective, initial data satisfying \eqref{3} is shown to also satisfy \eqref{4}, leading to finite-time blow-up. Therefore, both criteria for global existence and blow-up of solutions are equivalent.
	\end{proof}
	
	\paragraph{Acknowledgements.}
Li Chen is partially supported by the Deutsche Forschungsgemein-schaft (DFG, German Research Foundation)-547277619.
Yutian Lei is supported by the Natural Science Foundation of Jiangsu Province (BK20241878).

		{\sc Tiantian Zhou}
	
	Ministry of Education Key Laboratory for NSLSCS

	School of Mathematical Sciences
	
	Nanjing Normal University,
	Nanjing, 210023, China
	
	Email:zhoutiantiannj@163.com
	
	\vskip 5mm
	
	{\sc Li Chen}
	
	School of Business Informatics and Mathematics
	
	Universit\"{a}t Mannheim, 68131, Mannheim, Germany
	
	Email:li.chen@uni-mannheim.de
	
	\vskip 5mm
	
	{\sc Yutian Lei}
	
	Ministry of Education Key Laboratory for NSLSCS

	School of Mathematical Sciences
	
	Nanjing Normal University,
	Nanjing, 210023, China
	
	Email:leiyutian@njnu.edu.cn

\end{document}